\documentclass[11pt,a4paper]{article}
\usepackage{amsfonts}
\usepackage{epsfig}
\usepackage{amsmath}
\usepackage{amssymb}
\usepackage{color}
\usepackage{cite}

\usepackage[dvipdfm,CJKbookmarks,bookmarksopen=true,colorlinks=true,
linkcolor=blue,citecolor=blue,pdfstartview=FitH,pdftitle=title,pdfauthor=lixx]{hyperref}


\begin{document}%
\title{A Runge-Kutta discontinuous Galerkin scheme for hyperbolic conservation laws with discontinuous fluxes}
\author{Dian-liang Qiao$^{1,2}$
\quad Peng Zhang$^{1,5}$\thanks{Corresponding author. E-mail: pzhang@mail.shu.edu.cn.}\quad Zhi-Yang Lin$^{1}$\quad S.C. Wong$^3$\quad Keechoo Choi$^4$
\\
\small 1. Shanghai Institute of Applied Mathematics and Mechanics, Shanghai University,\\
\small Shanghai, P.R. China  \\
\small 2. Shanghai Industrial and Commercial Polytechnic, Shanghai, P.R. China,\\
\small 3. Department of Civil Engineering, The University of Hong Kong, Pokfulam Road,\\
\small Hong Kong SAR, P.R. China\\
\small 4. Department of Transportation Engineering, TOD-based Sustainable Urban Transportation Center,\\
\small Ajou University, Korea\\
\small 5. Shanghai Key Laboratory of Mechanics in
Energy Engineering
}
\date{ }%

%Create title of papers.
\maketitle

%\begin{abstract}

\noindent \textbf{Abstract}: The paper proposes a scheme by combining the Runge-Kutta discontinuous Galerkin method with a $\delta$-mapping algorithm for solving hyperbolic conservation laws with discontinuous fluxes. This hybrid scheme is particularly applied to nonlinear elasticity in heterogeneous media and multi-class traffic flow with inhomogeneous road conditions. Numerical examples indicate the scheme's efficiency in resolving complex waves of the two systems. Moreover, the discussion implies that the so-called $\delta$-mapping algorithm can also be combined with any other classical methods for solving similar problems in general.
\\

\noindent{{\bf Keywords}: $\delta$-mapping algorithm; elastic waves; multi-class traffic flow; Riemann problem; wave breaking }\\

\section{Introduction}

The standard hyperbolic conservation laws can be generally written in the following form \cite{Toro:1999,Shu:1998,Cockburn:2001}:
\begin{equation}\label{eq:1}
  u_{t}+f(u)_{x}=0,
\end{equation}
where $u=u(x,t)$ is an unknown variable or vector for solution, and $f(u)$ is the flux. However, considerable applications involve spatially varying fluxes, e.g., in flow through porous media, water wave equations, elastic waves in heterogeneous media, and traffic flow on an inhomogeneous road. See \cite{Lebacque:1996,LeVeque:2002a,LeVeque:2002b,LeVeque:2003,Bale:2003,Zhang:2003,Jin:2003,Jin:2009,Xu:2007,
Burger:2010a,Wang:2011,Wiens:2013} and the references therein for discussions of the problem. In this case, the equation or system for conservation is written as
    \begin{equation}\label{eq:2}
     u_{t}+f(u,\theta(x))_{x}=0,
     \end{equation}
where $\theta(x)$ is a known scalar or vector denoting some spatially varying parameters.

For standard conservation laws of Eq. (\ref{eq:1}), study of numerical schemes focuses on the capture of shocks. Although the first-order monotone scheme is able to resolve a shock, the profile can be over smoothed by numerical diffusions. The Godunov theorem suggests that a high-order accurate linear scheme can considerably reduce these diffusions. However, the dispersion that is due to the linearity yields spurious oscillations in the vicinity of a shock. Thus, nonlinearity was introduced in the high-order accurate scheme to suppress the oscillations, with the proposition of the total variation diminishing (TVD) scheme, the Runge-Kutta discontinuous Galerkin (RKDG) scheme, the weighted essentially non-oscillatory (WENO) scheme, and et al. See \cite{Bassi:1997,Shu:1998,Toro:1999,Cockburn:2001,Gottlieb:2001,LeVeque:2002a,Kubatko:2007} and the references therein for detailed discussions of the theory.

For conservation laws of Eq. (\ref{eq:2}), the aforementioned higher-order nonlinear schemes can be exploited for the capture of shocks. A straightforward treatment is to regard $\theta(x)$ as being continuous at the cell boundary $x_{j+1/2}$, or view $\theta(x)=\theta(x_{j+1/2})$ as being locally constant around $x_{j+1/2}$, and then directly applies these schemes by taking the numerical flux as $\hat{f}_{j+1/2}=\tilde{f}(u_{j+1/2}^-,u_{j+1/2}^+,\theta(x_{j+1/2}))$, where $\tilde{f}$ is a classical Riemann solver. However, such a treatment was indicated not to be consistent with the steady-state solution or stationary shock of Eq. (\ref{eq:2}), and oscillations were observed with relatively sharp change in $\theta(x)$ \cite{Zhang:2003,Zhang:2005a,Zhang:2005b}. We note that Eq. (\ref{eq:2}) usually gives a nonconstant steady-state solution $u=u(x)$, other than a trivial or constant solution that is implied in Eq. (\ref{eq:1}) or by setting $\theta(x)$ as being constant in Eq. (\ref{eq:2}).

Zhang and Liu \cite{Zhang:2003,Zhang:2005a} proposed a so-called $\delta$-mapping algorithm based on a thorough study of the characteristic theory under the scalar form of Eq. (\ref{eq:2}). The algorithm first assumes an intermediate state/value $\theta_{j+1/2}=\bar{\theta}(\theta_j,\theta_{j+1})$ of $\theta(x)$, which is somehow between the $j$-th and $(j+1)$-th cells, and then maps $u_j$ and $u_{j+1}$ onto the intermediate state $\theta_{j+1/2}$. The mapping is based on the fact that the flow $f(u,\theta(x))$ (other than the solution variable $u$) is constant in a characteristic. With the mapped values $\delta_{j+1/2}u_j$ and $\delta_{j+1/2}u_{j+1}$, the two adjacent solution states are \textquotedblleft unified" at a frozen state $\theta_{j+1/2}$, and a classical Riemann solver $\tilde{f}(\delta_{j+1/2}u_j,\delta_{j+1/2}u_{j+1},\theta_{j+1/2})$, e.g., the well known Godonov, Lax-Fridrichs, or Engquist-Osher flux is used to approximate the flux $f(u(x_{j+1/2},t),\theta(x_{j+1/2}))$ at the cell boundary. Since the Riemann solver depends on the two cell states $\theta_j$ and $\theta_{j+1}$ other than $\theta(x_{j+1/2})$, it is implied that the flux $f(u,\theta)$ is essentially discontinuous with respect to $\theta$ or $x$.

The so called $\delta$-mapping algorithm was further developed for solving Eq. (\ref{eq:2}), through combination with the RKDG scheme for the LWR model of traffic flow, with the WENO scheme for the elastic wave in heterogeneous media \cite{Xu:2007} and the multi-class model of traffic flow \cite{Zhang:2008}. These \textquotedblleft hybrid" schemes are different from the aforementioned \textquotedblleft straightforward treatment" in that $\delta_{j+1/2}u_i$ (other than $u_i$) were adopted in a classical numerical flux $\tilde{f}$, where $i$ refers to all involved cells for approximating $f(u,\theta)$ at $x=x_{j+1/2}$. These schemes were verified to be consistent with the stead-state flow or stationary shock of Eq. (\ref{eq:2}). We mention that other schemes, e.g., those developed for Eq. (\ref{eq:2}) in \cite{Lebacque:1996,Jin:2003,Jin:2009,Wang:2011,LeVeque:2002b,Bale:2003,Burger:2010a,Wiens:2013}, possess the same consistency and thus are able to well resolve the solution profiles despite much differences between their formulations and those in \cite{Zhang:2003,Zhang:2005a,Xu:2007,Zhang:2008}.

The present paper proposes a hybrid scheme for solving Eq. (\ref{eq:2}) by combining the $\delta$-mapping algorithm with the higher-order accurate RKDG scheme. Since the RKDG scheme (as well as TVD scheme) adopts a limiter that suggests nonlinearity or viscosity wherever near a shock, the $\delta$-mapping is also adopted in the limiter to maintain the aforementioned consistency with steady-state solutions or stationary shocks. Precisely, $u_{j\mp1}$ are replaced by $\delta_ju_{j\mp1}$ in the limiter referring to the $j$-th cell, where $\delta_j$ corresponds to the $\theta_j$ state. Although the discussion succeeds to that in \cite{Zhang:2005b}, we deal with the system more than the scalar equation, and focus on the multi-class traffic flow \cite{Wong:2002,ZhangM:2003,Zhang:2006,Donat:2008,Donat:2010,Ngoduy:2010,Ngoduy:2011,Basson:2009,Burger:2010a,
Burger:2010b,Burger:2011,Burger:2013,Chen:2012,Van:2013} and nonlinear elasticity in heterogeneous media \cite{LeVeque:2002a,LeVeque:2002b,LeVeque:2003,Bale:2003,Xu:2007}. The numerical results demonstrate that the scheme is robust in resolving the complex waves in the aforementioned problems, which are comparable with those given by the hybrid scheme that combines $\delta$-mapping and the fifth-order accurate WENO scheme in \cite{Xu:2007,Zhang:2008}, and those in \cite{LeVeque:2002b}.

The remainder of this paper is organized as follows. In Section 2, the RKDG scheme together with its combination with the $\delta$-mapping algorithm for the system of (\ref{eq:2}) is discussed in general. In Section 3, the aforementioned hybrid scheme is implemented with detailed discussions for elastic waves in heterogeneous media (Section 3.1) and for multi-class traffic flow (Section 3.2), respectively; numerical examples are presented in this section. We conclude the paper by Section 4.

\section{RKDG method combined with $\delta$-mapping}

\subsection{General account of DG space discretization}
A finite computational interval $[0,L]$ is uniformly divided into cells: $I_{j}=(x_{j-1/2},x_{j+1/2})$, with $\Delta_{j}=x_{j+1/2}-x_{j-1/2}$, and $x_{j}=(x_{j-1/2}+x_{j+1/2})/2$, $j=1,...,N$, which is shown by Fig. 1. For Eq. (\ref{eq:2}) with the initial condition:
\begin{equation}\label{eq:3}
 u(x,0)=u_{0}(x),
\end{equation}
we proceed the following. We multiply Eq. (\ref{eq:2}) with a test function $\omega(x)$, and integrate the resultant equation over $I_{j}$, which gives
\begin{equation*}
\int_{I_{j}}u_{t}\omega(x)dx+\int_{I_{j}}f_{x}(u,\theta)\omega(x)dx=0.
\end{equation*}
Then, we apply the integration by parts to the second term, and have
\begin{equation}\label{eq:4}
\int_{I_{j}}u_{t}\omega(x)dx-\int_{I_{j}}f(u,\theta)\omega_{x}(x)dx+f(u,\theta)\omega(x)|^{x_{j+1/2}}_{x_{j-1/2}}=0.
\end{equation}
Similar procedures are applied to Eq. (\ref{eq:3}), which yields
\begin{equation}\label{eq:5}
\int_{I_{j}}u(x,0)\omega(x)dx=\int_{I_{j}}u_{0}(x)\omega(x)dx.
\end{equation}

\begin{center}
 \epsfig{figure=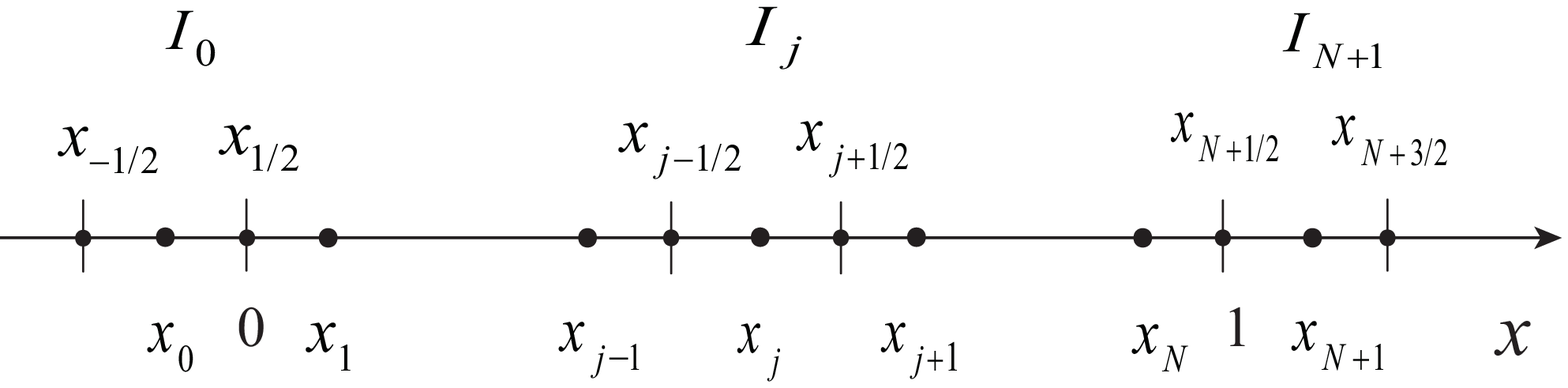,width=4.0 in}  \hspace{0.0 cm}\\
 \vspace{0.3cm}
 {\footnotesize Fig. 1 Cell division for space discretization}
\end{center}

Equations (\ref{eq:4}) and (\ref{eq:5}) are the weak formulations of Eqs. (\ref{eq:2}) and (\ref{eq:3}), on which the RKDG method is based for spatial discretization. Assume that $u_{h}(x,t)$ is an approximation to $u(x,t)$, where $u_{h}(x,t)\in P^{k}(I)$, and $P^{k}(I)$ is the space of piecewise polynomials of degree at most $k$. Then, in each cell $I_{j}$, $u_{h}(x,t)$ can be expressed as a linear combination:
\begin{equation*}
u_{h}(x,t)|_{I_{j}}=\sum_{l=0}^{k}u_{j}^{l}\varphi_{j}^{l}(x),
\end{equation*}
where $\{\varphi_{j}^{l}(x)\}_{l=0}^{k}$ is a set of bases of $P^{k}(I_{j})$. Usually, $\{\varphi_{j}^{l}(x)\}_{l=0}^{k}$ are taken as being orthogonal to each other under the $L^{2}$-norm, i.e.,
\begin{equation*}
 \varphi_{j}^{l}(x)=L_{l}(\frac{2(x-x_{j})}{\Delta_{j}}), L_{l}(s)=\frac{1}{2^{l}l!}\frac{d^{l}}{ds^{l}}[(s^{2}-1)^{l}],
\end{equation*}
where $L_{l}$ represent the $Legendre$ polynomials. By taking $\omega(x)=\varphi_{j}^{l}(x)$, and replacing $u(x,t)$ with $u_{h}(x,t)$, Eq. (\ref{eq:4}) leads to the following ordinary differential equations (ODEs):
\begin{equation}\label{eq:6}
\frac{d}{dt}u_{j}^{l}(t)=\frac{2l+1}{\Delta_{j}}(\int_{I_{j}}f(u_{h}(x,t),\theta(x))(\varphi_{j}^{l}(x))_{x}dx
-\hat{f}_{j+1/2}+(-1)^{l}\hat{f}_{j-1/2}),
\end{equation}
for solving the coefficients $u_{j}^{l}(t)$. Here, the numerical flux $\hat{f}_{j\pm1/2}$ is used to approximate the flow $f(u,\theta)$ at the cell boundary $x=x_{j\pm1/2}$, through the replacement of $f(u_{h}(x_{j\pm1/2},t),\theta(x_{j\pm1/2}))$ by $\hat{f}_{j\pm1/2}$. Similar operations are implemented on (\ref{eq:5}), which help derive the initial values of these coefficients with
 \begin{equation}\label{eq:7}
 u_{j}^{l}(0)=\frac{2l+1}{\Delta_{j}}\int_{I_{j}}u_{0}(x)\varphi_{j}^{l}(x)dx.
   \end{equation}

The approximation is theoretically of the $(k+1)$-th order of accuracy. To ensure the same order of accuracy, all integrals are computed by the Gauss formula with sufficiently high accuracy (e.g., the two-point formula for $k=1$). We refer the reader to \cite{Bassi:1997,Cockburn:2001,Gottlieb:2001,Kubatko:2007,Qiao:2014} for general account of the formulation.

\subsection{Derivation of numerical fluxes}
We omit the subscript \textquotedblleft $h$" in $u_h$, and denote by $u^{\pm}_{j+1/2}(t)=u_{h}^{\pm}(x_{j+1/2},t)$, for the discussion in this section. For $\theta$ being considered continuous at the cell boundary $x=x_{j+1/2}$, the flow $f(u_{h}(x_{j\pm1/2},t),\theta(x_{j\pm1/2}))$ could be approximated by an classical numerical flux, which takes the following form:
    \begin{equation*}
    \hat{f}_{j+1/2}=\tilde{f}(u_{j+1/2}^{-}(t),u_{j+1/2}^{+}(t),\theta(x_{j+1/2})),
    \end{equation*}
and which refers to the aforementioned \textquotedblleft straightforward treatment". However, this approximation suggests non-consistency with the steady-state or stationary shock solution of Eq. (\ref{eq:2}), and non-physical oscillations for sharp change in $\theta(x)$ \cite{Zhang:2003,Zhang:2005a,Zhang:2005b}. Therefore, $\theta(x)$ should be viewed as being discontinuous at $x=x_{j+1/2}$, and the numerical flux should take the following form:
    \begin{equation}\label{eq:8}
      \hat{f}_{j+1/2}=\hat{f}(u_{j+1/2}^{-},\theta^{-}(x_{j+1/2});u_{j+1/2}^{+},\theta^{+}(x_{j+1/2})) \end{equation}
For the studied problems in this paper, $\theta(x)$ is piece-wise constant, thus we simply set $\theta^{-}(x_{j+1/2})=\theta(x_{j})\equiv\theta_j$, and $\theta^{+}(x_{j+1/2})=\theta(x_{j+1})\equiv\theta_{j+1}$. In an otherwise case, $\theta(x)$ could be properly approximated by a polynomial.

To alternatively exploit a classical numerical flux $\tilde{f}$ for the definition of $\hat{f}$ in Eq. (\ref{eq:8}), we proceed the following. We choose a certain intermediate state $\theta_{j+1/2}$ between $\theta_j$ and $\theta_{j+1}$, $\theta_{j+1/2}=\bar{\theta}(\theta_j,\theta_{j+1})$, where $\bar{\theta}$ is an average between $\theta_j$ and $\theta_{j+1}$, such that $\bar{\theta}(\theta,\theta)=\theta$. Then, $u^{-}_{j+1/2}$ and $u^{+}_{j+1/2}$ are \textquotedblleft unified" by mapping them onto $\theta_{j+1/2}$ state with the mapped values $\delta_{j+1/2}u^{\mp}_{j+1/2}$, and $\hat{f}$ is given by
      \begin{equation}\label{eq:9}
      \hat{f}(u_{j+1/2}^{-},\theta_{j};u_{j+1/2}^{+},\theta_{j+1})
      =\tilde{f}(\delta_{j+1/2}u^{-}_{j+1/2},\delta_{j+1/2}u^{+}_{j+1/2},\theta_{j+1/2}).
      \end{equation}
The mapping $\delta_{j+1/2}$ is defined by the following. Given $u^{\mp}_{j+1/2}$, we find $\delta_{j+1/2}u^{\mp}_{j+1/2}$, which maximizes $\gamma\in(-\infty, 1]$ under the restrictions:
\begin{equation}\label{eq:10}
f(\delta_{j+1/2}u^{\mp}_{j+1/2},\theta_{j+1/2})=\gamma f(u^{\mp}_{j+1/2},\theta_{j+1/2\mp1/2}),
\end{equation}
and
  \begin{equation}\label{eq:11}
   \lambda_{l}(\delta_{j+1/2}u^{\mp}_{j+1/2},\theta_{j+1/2})\lambda_{l}(u^{\mp}_{j+1/2},\theta_{j+1/2\mp1/2})\geq0,
\end{equation}
where $\{\lambda_l\}$ are the eigenvalues of Eq. (\ref{eq:2}) with $\theta$ being fixed; moreover, we set
\begin{align}\label{eq:12}
\left \{\aligned
 &\lambda_{l}(\delta_{j+1/2}u^{-}_{j+1/2},\theta_{j+1/2})\geq0,\;\;\; if\; \lambda_{l}(u^{-}_{j+1/2},\theta_{j})=0,\\
 &\lambda_{l}(\delta_{j+1/2}u^{+}_{j+1/2},\theta_{j+1/2})\leq0,\;\;\; if\; \lambda_{l}(u^{+}_{j+1/2},\theta_{j+1})=0.
 \endaligned
 \right.
 \end{align}

By Eq. (\ref{eq:10}), we attempt to equalize the two flows with $\gamma=1$, or at least maximize the flow at $\theta_{j+1/2}$. This is in accordance with the \textquotedblleft supply-demand" concept used in the theory of fluid dynamics or traffic flow \cite{Lebacque:1996,Jin:2009}. For the first equation of (\ref{eq:10}), the demand $f(u^{-}_{j+1/2},\theta_j)$ is fully satisfied at $\theta_{j+1/2}$ with $\gamma=1$, if it (or its components) is not larger than the capacity (or capacities) of $f$ at $\theta_{j+1/2}$; otherwise, it is partly satisfied by reaching the capacity (or capacities) at $\theta_{j+1/2}$. For the second equation of (\ref{eq:10}), the supply $f(u^{+}_{j+1/2},\theta_{j+1})$ is fully available at $\theta_{j+1/2}$ with $\gamma=1$, if it (or its components) does not exceed the capacity (or capacities) of $f$ at $\theta_{j+1/2}$; otherwise, it is partly available by reaching the capacity (or capacities) at $\theta_{j+1/2}$. Eqs. (\ref{eq:11})-(\ref{eq:12}) imply a so called wave entropy condition, i.e., valid information in a certain characteristic should go forward or backward without turning back. As will be shown in Sections 3 and 4, Eqs. (\ref{eq:10})-(\ref{eq:12}) are applicable to the discussed problems.

\subsection{Time discretization and limiter}
Equations (\ref{eq:6})-(\ref{eq:7}) can be rewritten as the following ODEs:
\begin{equation}\label{eq:13}
\frac{du_{h}}{dt}|_{I_{j}}=L_{h}(u_{h}), u_{h}|_{I_{j}}(0)=\sum_{l=0}^{k}u_{j}^{l}(0)\varphi_{j}^{l}(x),\;j=1,...,N,
\end{equation}
for which the $(k+1)$-th order TVD Runge-Kutta time discretization is adopted. Note that we retrieve the subscript \textquotedblleft $h$" in $u_h$. The procedure is briefed in the following.

For a division $\{t^{n}\}_{n=0}^{M}$ of the time interval $[0,T]$, where $t^{0}=0$, and $\Delta t^{n}=t^{n+1}-t^{n}$, we set ${u}_{j}^{0}=\Lambda\Pi_{h}^{k}{u_{h}}_{j}(0)$. Then , for $n=0,...,M-1$, ${u_{h}}_{j}^{n+1}$ are computed as follows:\\
(i) Set ${u_{h}}_{j}^{(0)}={u_{h}}_{j}^{n}$;\\
(ii) For $i=1,...,k+1$, compute the intermediate functions:
\begin{equation*}
{u_{h}}_{j}^{(i)}=\Lambda\Pi_{h}^{k}\{\sum_{l=0}^{i-1}\alpha_{il}{u_{h}}_{j}^{(l)}+\beta_{il}\Delta t^{n}L_{h}({u_{h}}_{j}^{(l)})\};
\end{equation*}
(iii) Set ${u_{h}}_{j}^{n+1}={u_{h}}_{j}^{(k+1)}$.

All the procedure (together with the parameters $\alpha_{il}$ and $\beta_{il}$) is the same as that in \cite{Cockburn:2001}, except that the slope limiter $\Lambda\Pi_{h}^{k}$ is redesigned to guarantee the scheme's consistency with a steady-state solution or stationary shock.

For $k=1$, the slope limiter that acts on the piecewise linear solution
   \begin{equation*}
   u_{h}|_{I_{j}}=\bar{u}_{j}+\frac{2u_{j}^{1} }{\Delta_{j}}(x-x_{j}),\;j=1,...,N,
    \end{equation*}
is defined by
\begin{equation}\label{eq:14}
\Lambda\Pi_{h}^{1}u_{h}|_{I_{j}}=\bar{u}_{j}+\frac{2}{\Delta_{j}}
m(u_{j}^{1},\delta_{j}\bar{u}_{j+1}-\bar{u}_{j},\bar{u}_{j}-\delta_{j}\bar{u}_{j-1})(x-x_{j}),\;j=1,...,N. \end{equation}
Here, the definition of $\delta_j$ is similar to that of $\delta_{j+1/2}$ by Eqs. (\ref{eq:9})-(\ref{eq:12}). However, $\delta_j$ is used to map the two adjacent averages $\bar{u}_{j\mp1}$ onto $\theta_j$ state for comparison with the average $\bar{u}_{j}$. By Eq. (\ref{eq:14}) the average of $u_{h}|_{I_{j}}$ remains the same for the purpose of conservation; however, the slope or change in $u_{h}|_{I_{j}}$ is possibly limited which with resultant artificial viscosities helps suppress non-physical oscillations. Eq. (\ref{eq:14}) gives the following cell boundary values,
\begin{equation}\label{eq:15}
u_{j+1/2}^{-}=\bar{u}_{j}+m(u_{j+1/2}^{-}-\bar{u}_{j},\delta_{j}
\bar{u}_{j+1}-\bar{u}_{j},\bar{u}_{j}-\delta_{j}\bar{u}_{j-1}),
\end{equation}
\begin{equation}\label{eq:16}
u_{j-1/2}^{+}=\bar{u}_{j}-m(\bar{u}_{j}-u_{j-1/2}^{+},\delta_{j}
\bar{u}_{j+1}-\bar{u}_{j},\bar{u}_{j}-\delta_{j}\bar{u}_{j-1}),
\end{equation}
which are used in Eqs. (\ref{eq:10})-(\ref{eq:12}). The minmod function $m$ is defined by
\begin{align*}
m(a_{1},a_{2},a_{3})=\left \{\aligned
 &s\mathop{\min}\limits_{1\leq n\leq3}|a_{n}| ,\quad if\;s=sign(a_{1})=sign(a_{2})=sign(a_{3}),\\
 &0,\quad otherwise.
 \endaligned
 \right.
 \end{align*}

For $k>1$, the limiter $\Lambda\Pi_{h}^{k}$ is defined based on the definition of $\Lambda\Pi_{h}^{1}$, which follows almost the same steps in \cite{Gottlieb:2001,Cockburn:2001} (see also \cite{Zhang:2005b}).
To ensure the numerical stability of a scheme, the time step should satisfy the following CFL condition\cite{Gottlieb:2001,Cockburn:2001}:
  \begin{equation}\label{eq:17}
  \Delta t^{(n)}\leq C\frac{\Delta_{j}}{\alpha^{(n)}},\;\;C=\frac{1}{2k+1},
  \end{equation}
where, $\alpha^{(n)}=\max_{j}\max_{i}\{|\lambda_{1}(\delta_{j}u_{i}^{(n)},\theta_{j})|,...,|\lambda_{r}(\delta_{j}u_{i}^{(n)},\theta_{j})|\}$.

The formulation is almost the same as that in the standard RKDG scheme \cite{Cockburn:2001}, except that $u_j$ and $u_{j\mp1}$ in the limiter are replaced by $\delta_{j}u_j$ and $\delta_{j}u_{j\mp1}$. It can be easily verified that all the procedures guarantee the scheme's consistency with the steady-state solution or stationary shock of system (\ref{eq:2}). On the other hand, a \textquotedblleft straightforward treat" that is without $\delta_{j}$ in Eqs. (\ref{eq:14})-(\ref{eq:17}) does not possess such a consistency and thus would suggest non-physical oscillations. We refer the reader to \cite{Zhang:2005b} for a similar discussion for the scalar equation of (\ref{eq:2}).

\section{Numerical implementation}

\subsection{Application to elastic wave equations}
We consider the following elastic wave equations:
\begin{align}\label{eq:18}
& \varepsilon(x,t)_{t}-v(x,t)_{x}=0,\\ \label{eq:19}
& (\rho(x)v(x,t))_{t}-\sigma(\varepsilon(x,t),K(x))_{x}=0,
\end{align}
where the strain $\varepsilon(x,t)$ and the velocity $v(x,t)$ are the unknowns, and the stress-strain relation is given by
   \begin{equation*}
   \sigma(\varepsilon,K)=K\varepsilon+\beta K^{2}\varepsilon^{2}.
     \end{equation*}
By setting $u=(\varepsilon,q)^{T}$, $\theta(x)=(\rho(x),K(x))^{T}$, $q=\rho v$, and  $f(u,\theta)=(-q/\rho,-\sigma(\varepsilon,K))^{T}$, the system of (\ref{eq:18})-(\ref{eq:19}) takes the form of Eq. (\ref{eq:2}). For $\theta(x)$ being fixed, the two eigenvalues of the system are easily shown as $\lambda_{1}=-c$, and $\lambda_{2}=c$, where $c=\sqrt{\sigma_{\varepsilon}/\rho}$ is the sonic speed.

To implement the numerical scheme in Section 2 for solving the system of (\ref{eq:18})-(\ref{eq:19}), we only need verify that Eqs. (\ref{eq:10})-(\ref{eq:12}) are applicable. We can actually choose $\gamma=1$, such that Eq. (\ref{eq:10}) is always solvable with
   \begin{equation*}
\delta_{j+1/2}q^{\mp}_{j+1/2}=\frac{\rho_{j+1/2}q^{\mp}_{j+1/2}}{\rho_{j+1/2\mp1/2}},\;\; \delta_{j+1/2}\varepsilon^{\mp}_{j+1/2}=\frac{(1+4\beta\sigma(\varepsilon^{\mp}_{j+1/2},K_{j+1/2\mp1/2}))^{1/2}-1}{2\beta K_{j+1/2}},
\end{equation*}
where the intermediate state $(\rho_{j+1/2},K_{j+1/2})=(\rho_{j+1},K_{j+1})$, and we set $\beta=0.3$ in the simulation. Eqs. (\ref{eq:11})-(\ref{eq:12}) are self-evident in that there always hold $\lambda_1<0$, and $\lambda_2>0$. The mapped values used in Eqs. (\ref{eq:14})-(\ref{eq:17}) can be similarly derived.

The system of (\ref{eq:18})-(\ref{eq:19}) is used to model a compound material consisting of alternating layers of two different materials, in which case the density $\rho(x)$ and the modulus $K(x)$ are taken as piecewise constants for $x\in[0,300]$, with
  \begin{align*}
(\rho(x),K(x))=\left \{\aligned
 &(\rho_{A},K_{A}),\;\;2k\leq x<2k+1,\\
 &(\rho_{B},K_{B}),\;\;2k+1\leq x<2k+2,
 \endaligned
 \right. \;k=0,\cdots,149.
 \end{align*}
Precisely, we set $(\rho_{A},K_{A})=(1,1)$, and $(\rho_{B},K_{B})=(3,3)$. Initially,
$\varepsilon(x,0)=0$, and $v(x,0)=0$; a perturbation is put on the left boundary with $v(0,t)=-0.2(1+\cos(\pi(t-30)/30))$, for $t\leq60$, which is cleared up with $v(0,t)=0$, for $60<t<70$. Thereafter, the periodic boundary conditions are applied to observe the development of the perturbation, which is actually regarded as going forward to the infinite.

\begin{center}
\scriptsize(a)\epsfig{figure=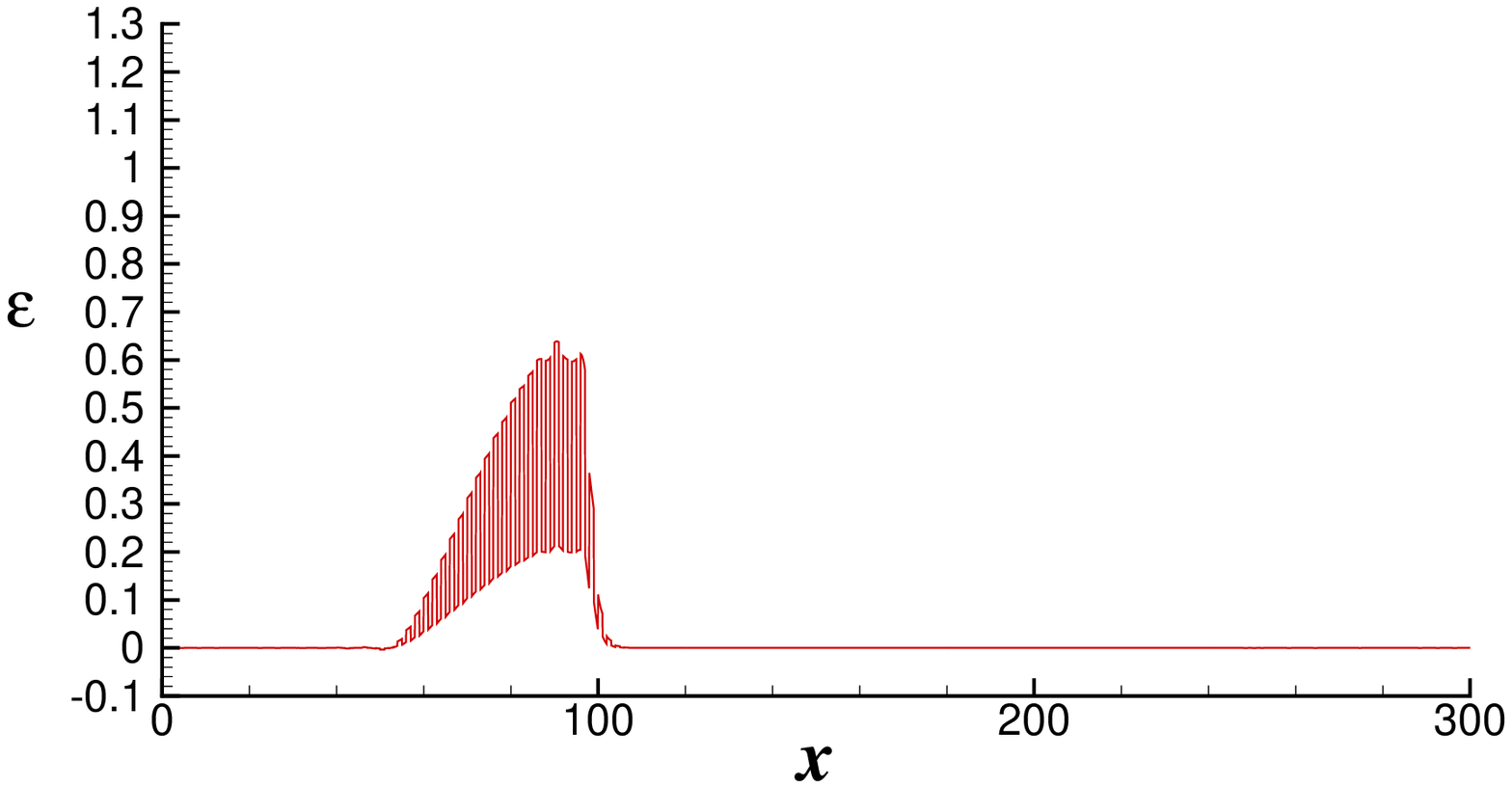,width=3.0 in}  \hspace{0.0 cm}
\scriptsize(b)\epsfig{figure=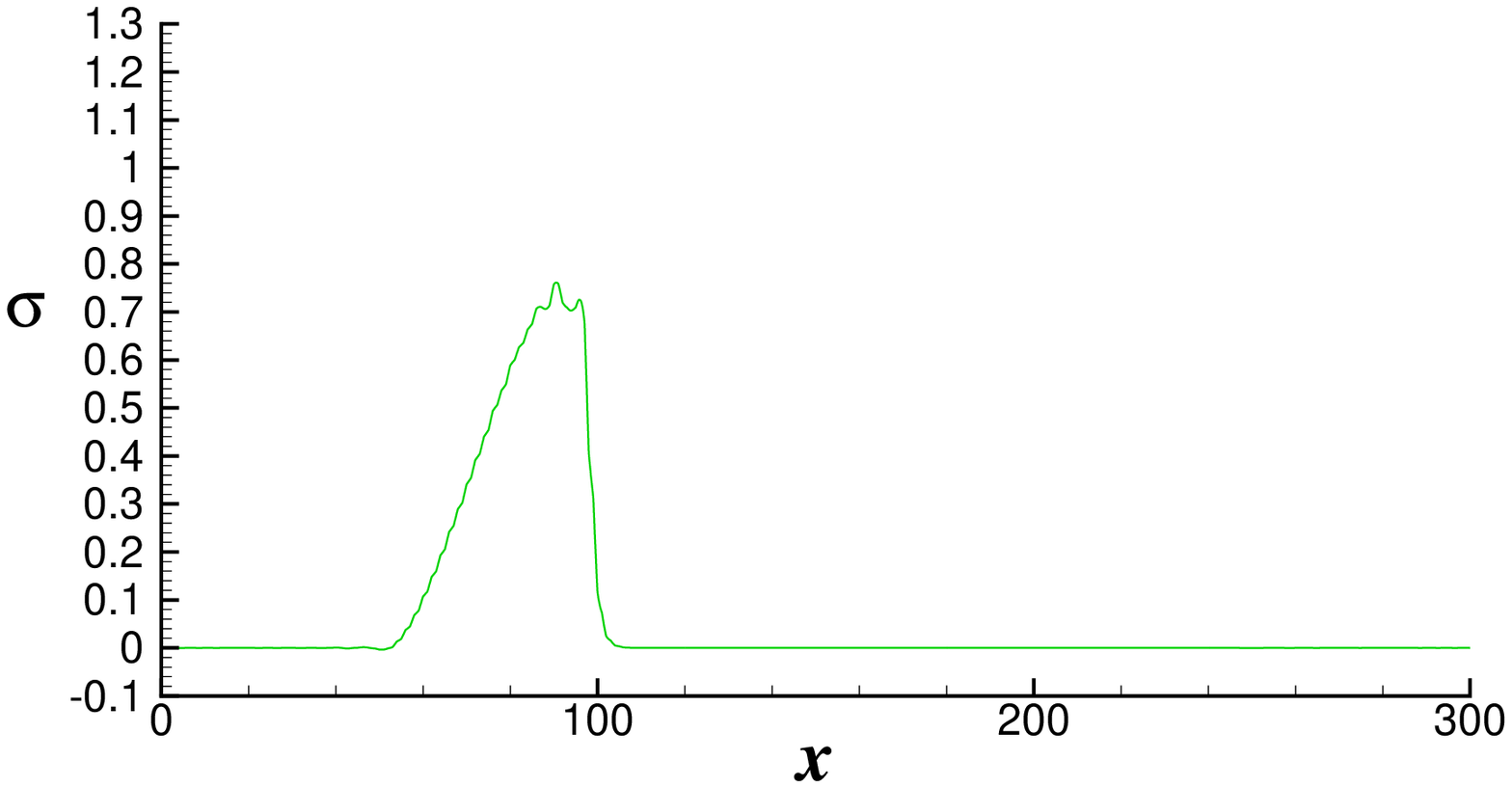,width=3.0 in}  \hspace{0.0 cm}
\scriptsize(c)\epsfig{figure=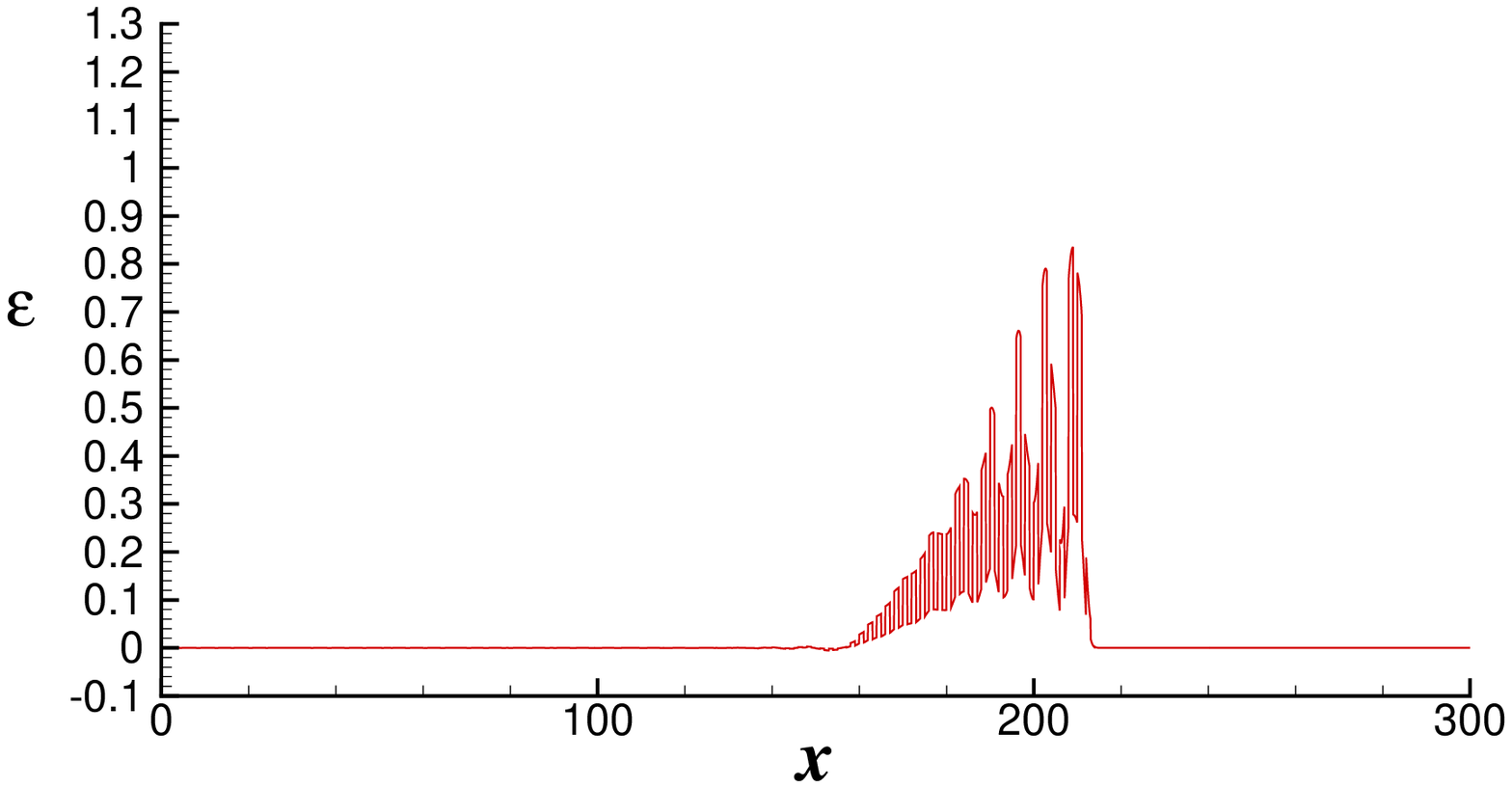,width=3.0 in}  \hspace{0.0 cm}
\scriptsize(d)\epsfig{figure=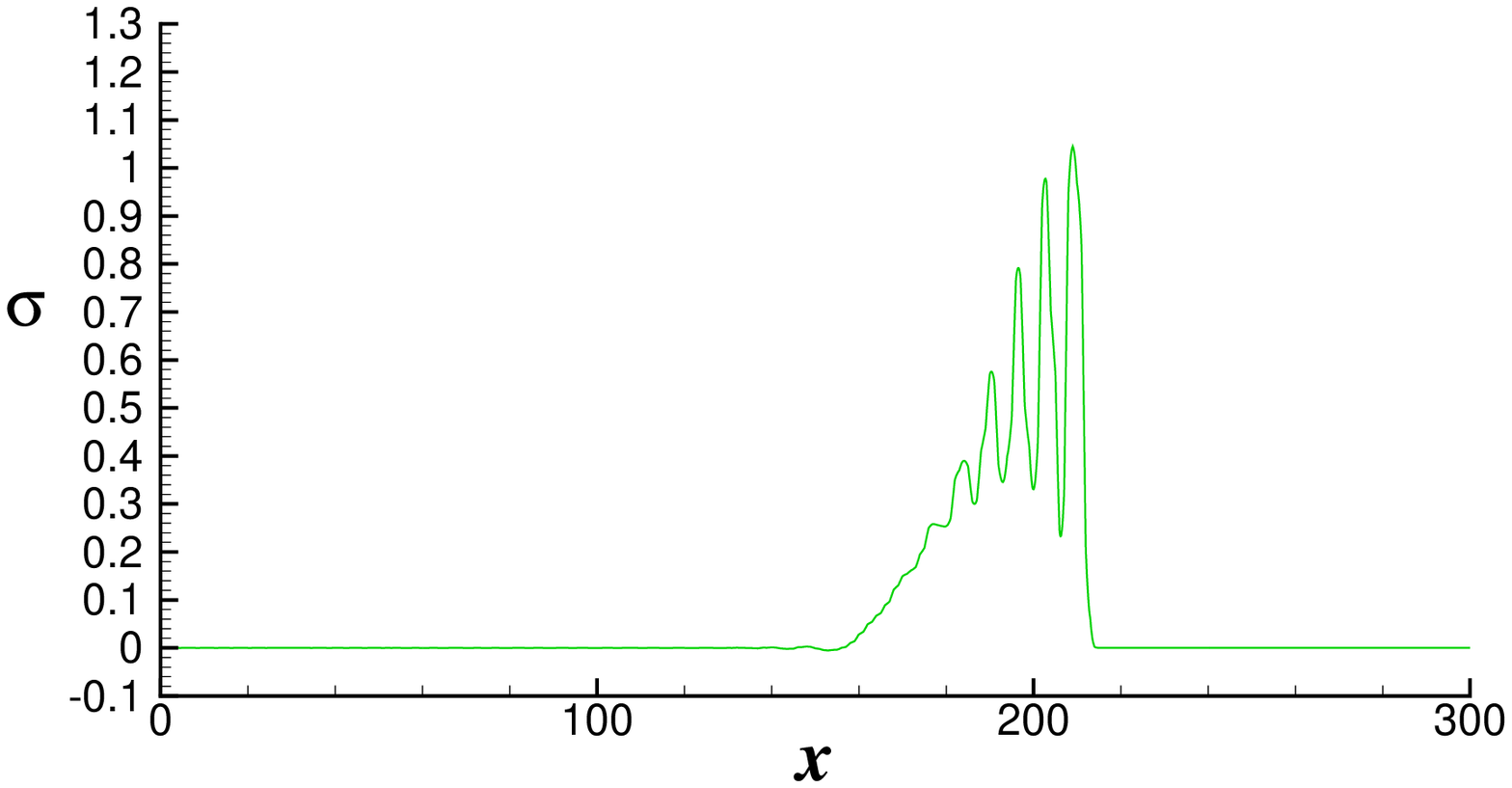,width=3.0 in}  \hspace{0.0 cm}\\
 \vspace{0.3cm}
 {\footnotesize Fig. 2 The strain and stress, (a)-(d) at $t=120$ and $t=240$.}
\end{center}

\begin{center}
\scriptsize(a)\epsfig{figure=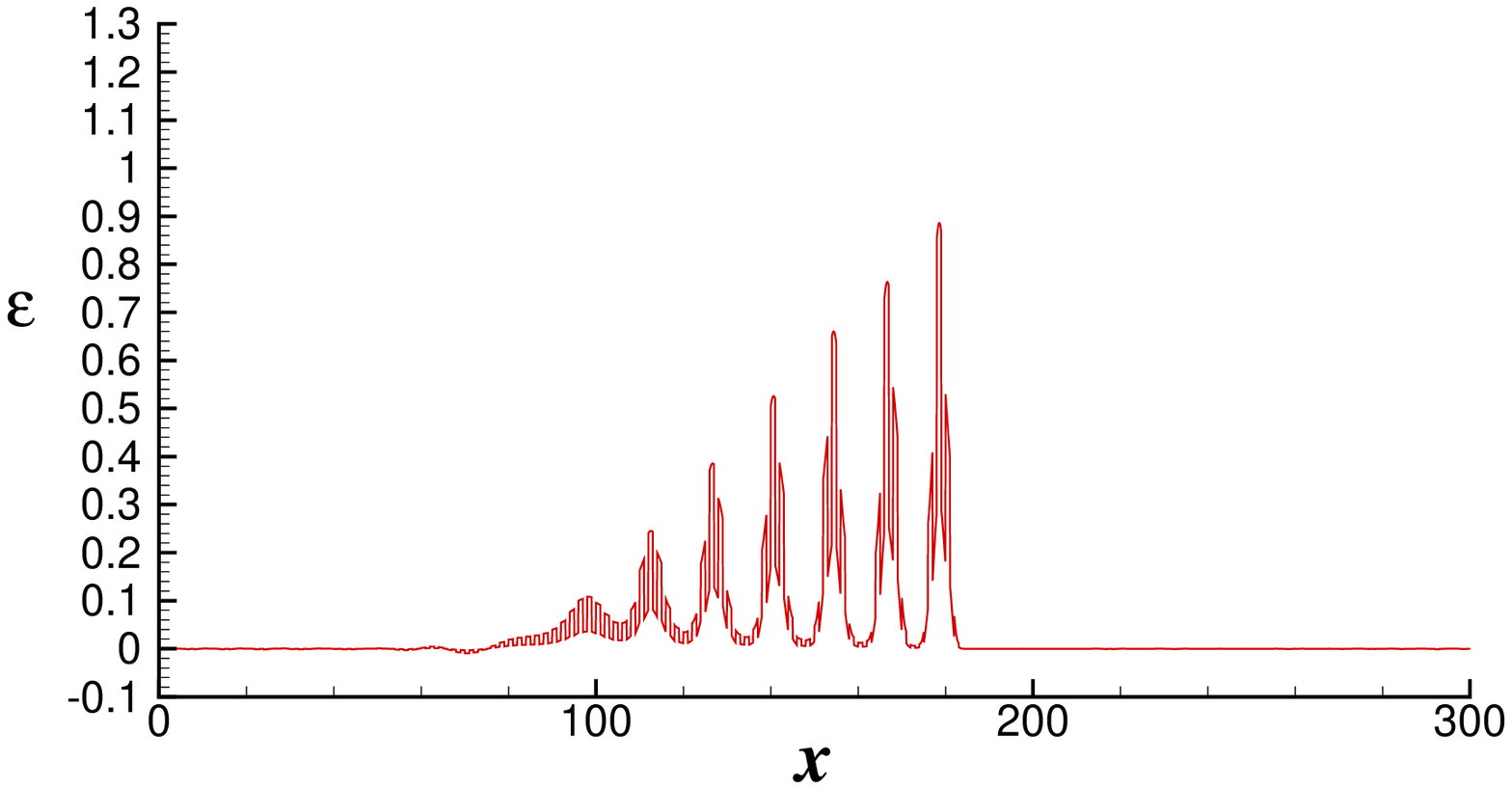,width=3.0 in}  \hspace{0.0 cm}
\scriptsize(b)\epsfig{figure=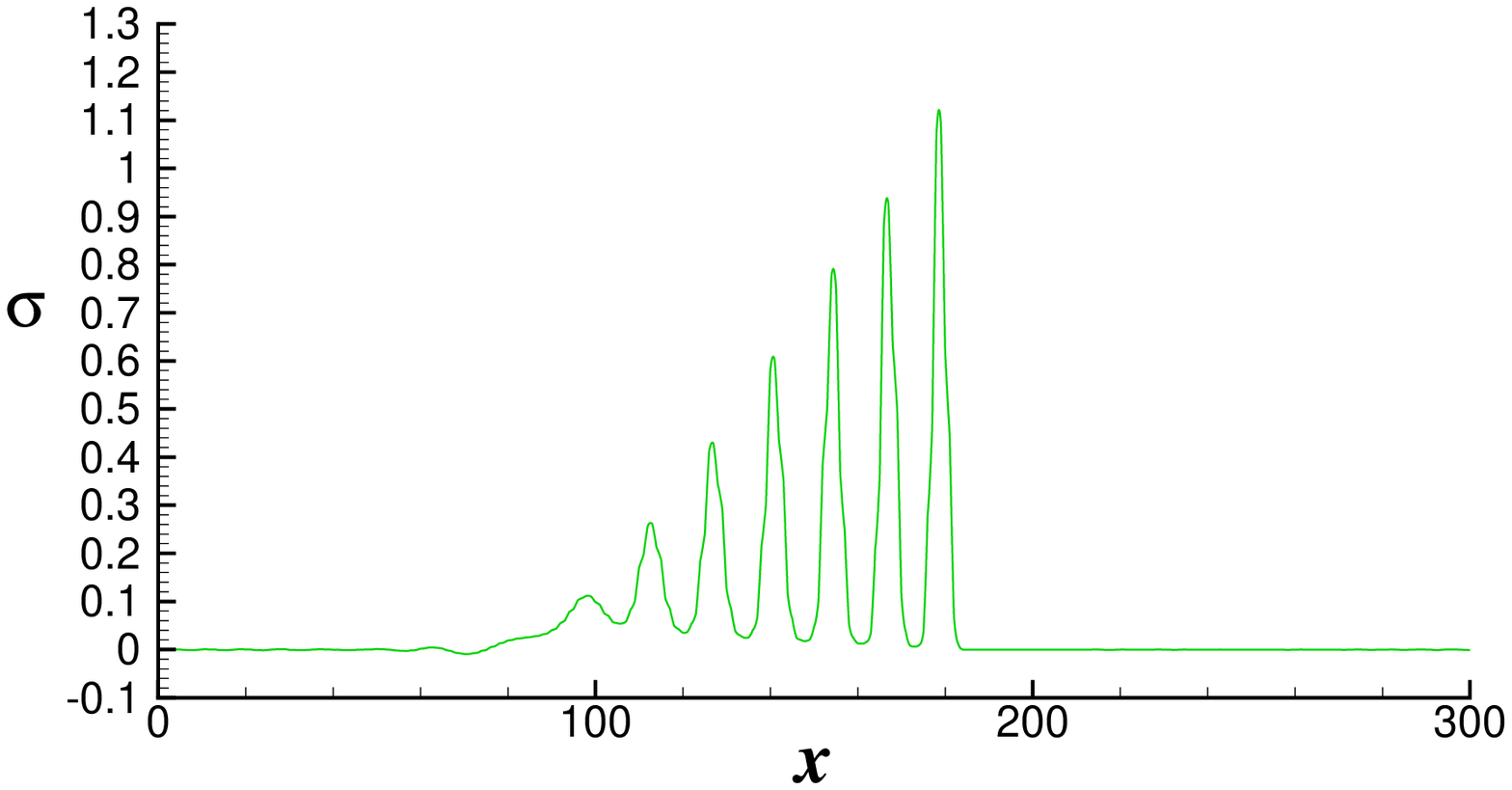,width=3.0 in}  \hspace{0.0 cm}
\scriptsize(c)\epsfig{figure=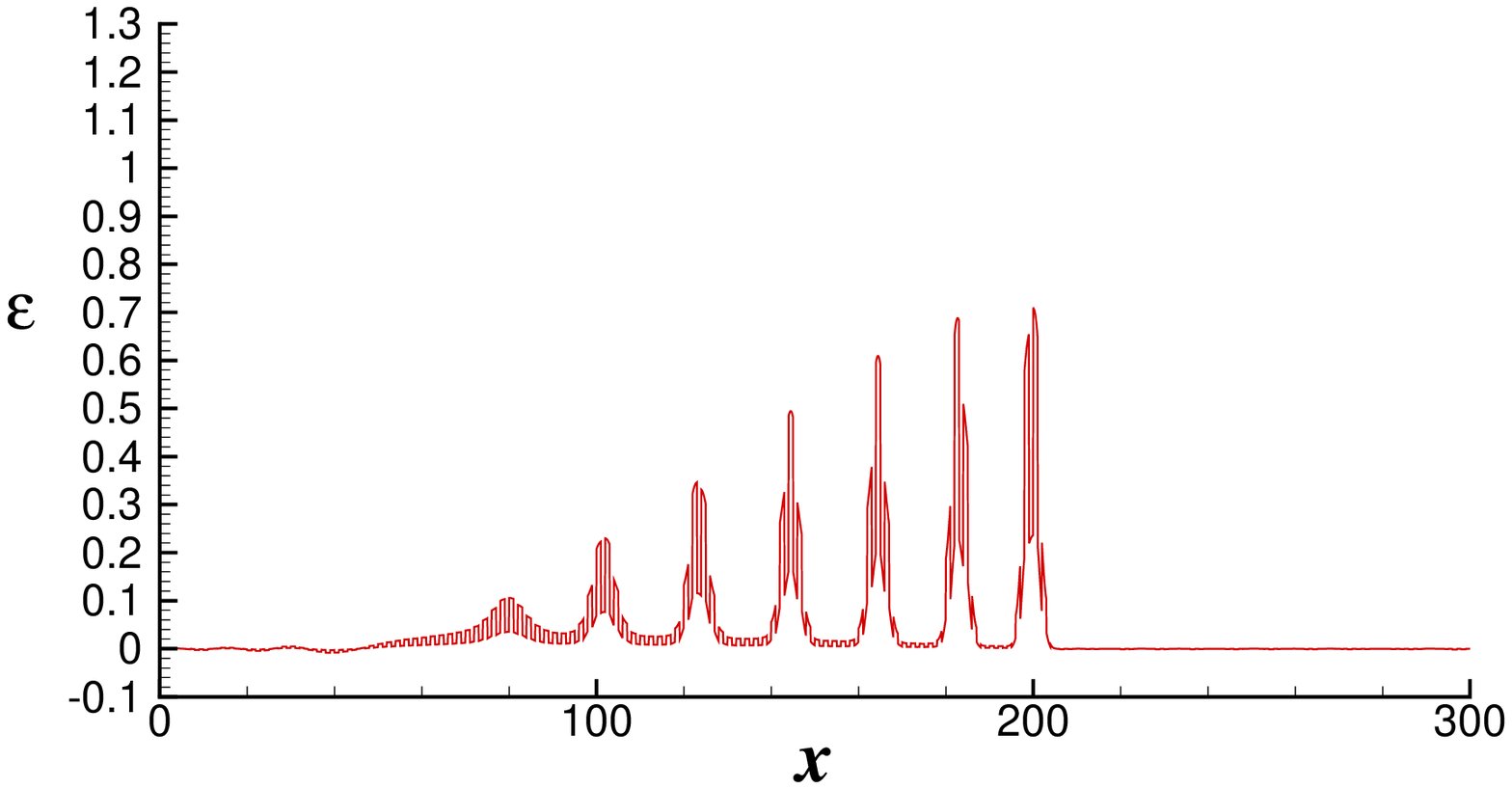,width=3.0 in}  \hspace{0.0 cm}
\scriptsize(d)\epsfig{figure=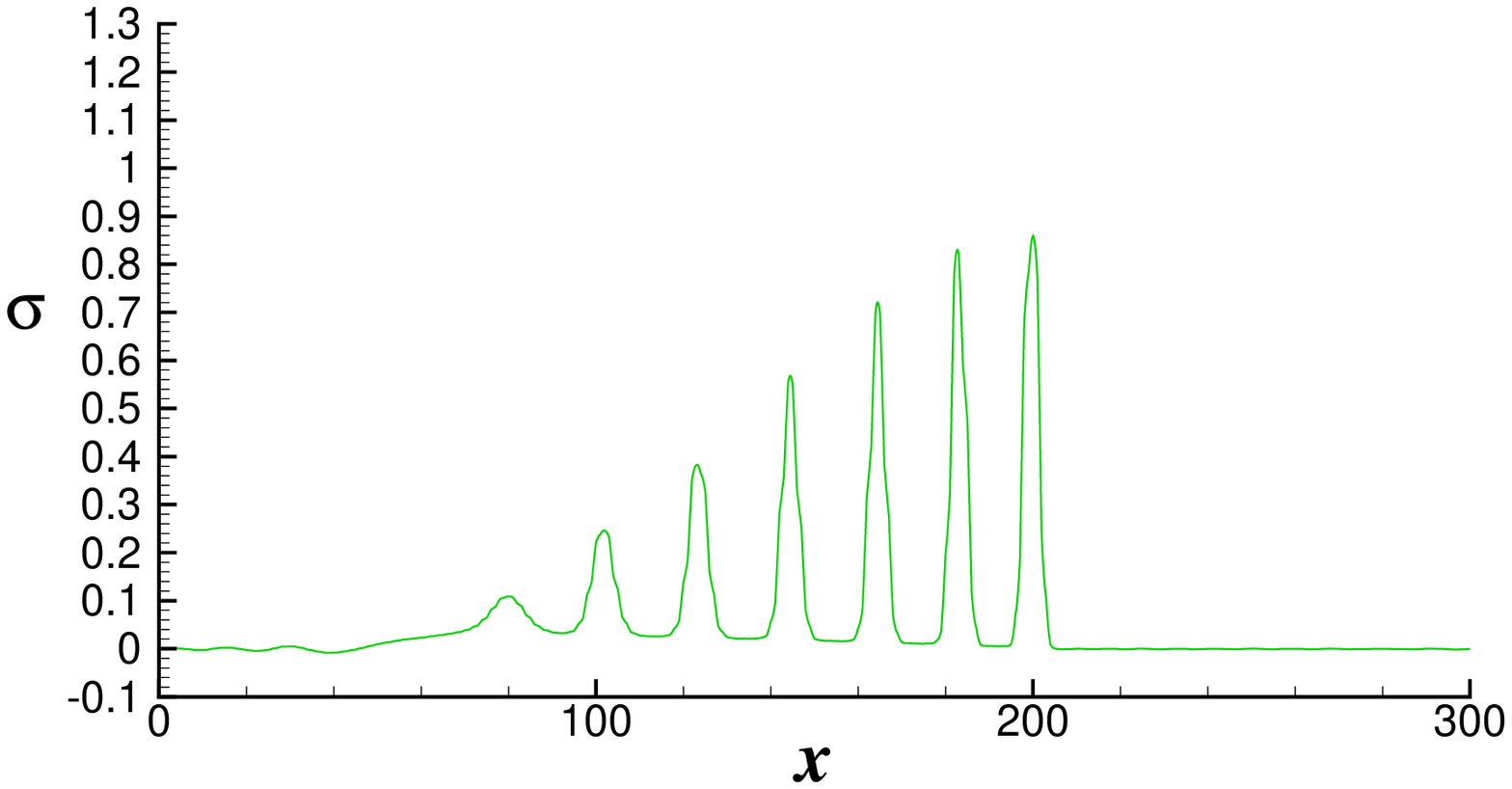,width=3.0 in}  \hspace{0.0 cm}
\scriptsize(e)\epsfig{figure=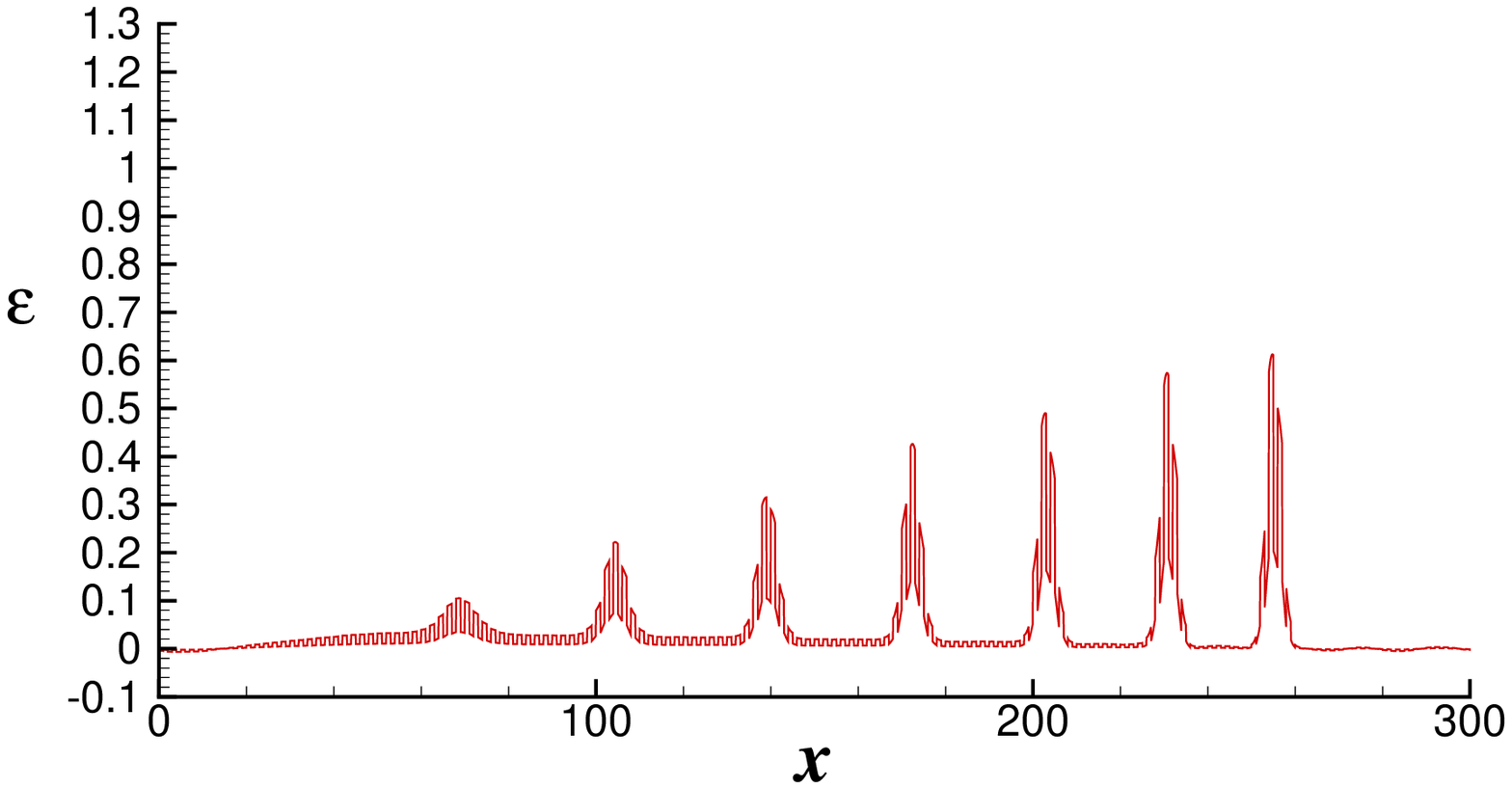,width=3.0 in}  \hspace{0.0 cm}
\scriptsize(f)\epsfig{figure=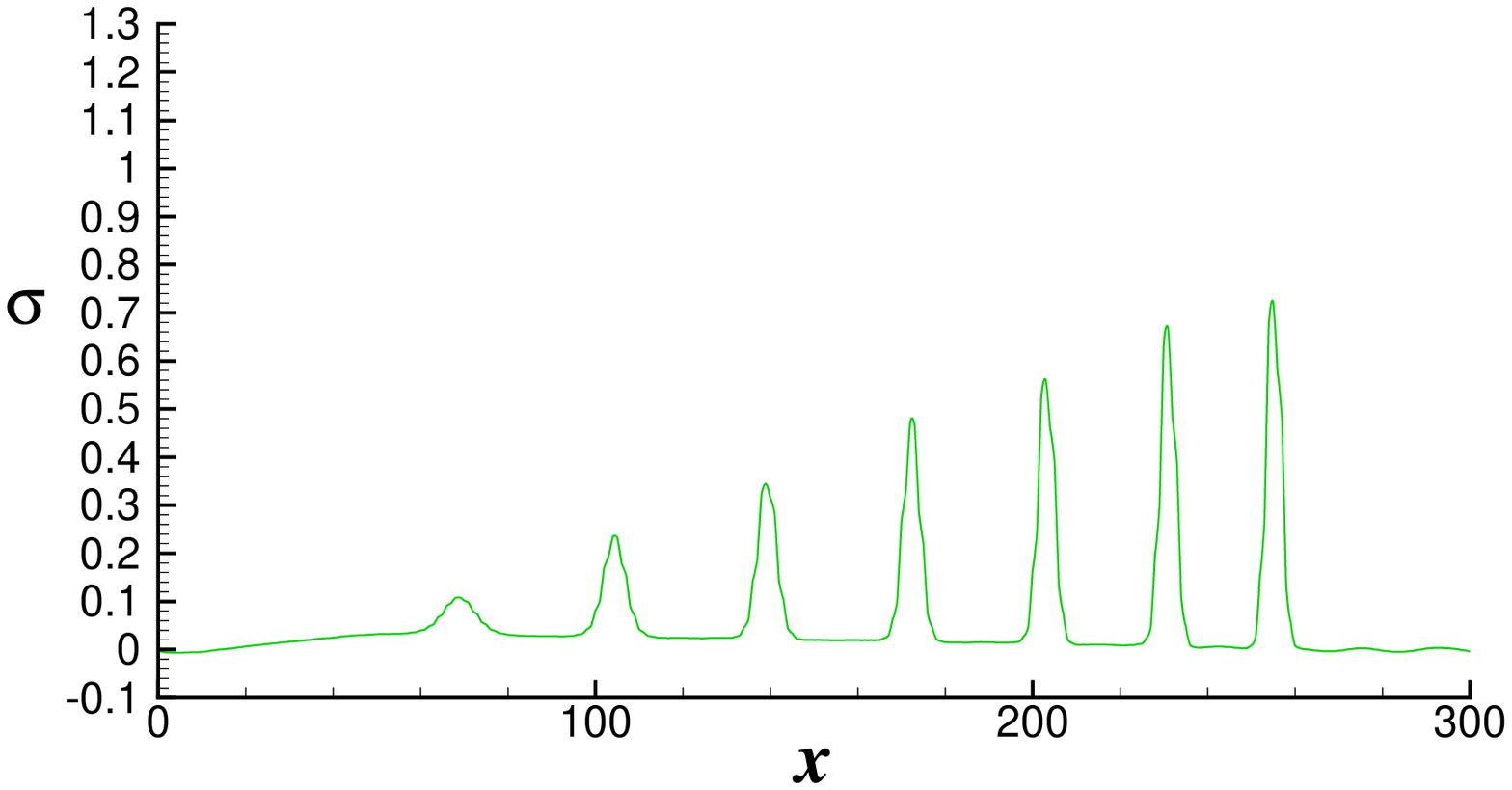,width=3.0 in}  \hspace{0.0 cm}\\
 \vspace{0.3cm}
 {\footnotesize Fig. 3 The strain and stress at (a), (b) $t=840$; (c), (d) $t=1500$; and (e), (f) $t=2850$.}
\end{center}

Figures 2 and 3 show that the perturbation breaks up into a series of solitary waves, which have similar shapes. However, the magnitudes and propagation speeds are different; the larger the magnitude, the faster the propagation. These results indicate the efficiency of the scheme in resolving nonlinear elastic wave in heterogeneous media. See also \cite{LeVeque:2002b,Xu:2007} for comparison.

\subsection{Application to multi-class traffic flow}
The multi-class model proposed in \cite{Wong:2002} was extended to describe traffic flow on an inhomogeneous highway road, which takes the same form of Eq. (\ref{eq:2}), i.e.,
   \begin{equation}\label{eq:20}
   u_t+f(u,\theta(x))_x=0,
   \end{equation}
with $u=(u_{1},...,u_{m})^{T}=(a\rho_{1},...,a\rho_{m})^{T}$, $f(u,\theta)=(f_{1}(u,\theta),...,f_{m}(u,\theta))^{T}$, and $\theta(x)=(a(x),b_{1}(x),...,b_{m}(x))^{T}$. Here, $a(x)$ is the number of lanes, $\rho_l$ and $f_l$ are the average density per lane and the flow of the $l$-th vehicular species, and $b_{l}(x)=v_{l,f}(x)/v_{f}$ is the scaled free flow velocity with $v_{f}=\max_x\max_{1\leq l\leq m}\{v_{l,f}(x)\}$.

The average velocity $v_l\equiv f_l/u_l$ of the $l$-th vehicular species is taken as $b_{l}(x)v(\rho)$, where $\rho=\sum_{l=1}^m\rho_l$, and $v(\rho)$ is a velocity-density relationship. Thus, we have
 \begin{equation*}
    v_l(x,t)=b_{l}(x)v(\rho),\;\;f_l(u,\theta)=b_{l}(x)u_lv(\rho)=b_{l}u_lv(\sum_{l=1}^mu_l/a).
    \end{equation*}
We adopt the same velocity-density relationship as that in \cite{Zhang:2008}:
   \begin{equation*}
   v(\rho)=v_f(1-\frac{\rho}{\rho_{jam}}),
   \end{equation*}
and assume that the $m$ vehicular species are divided, such that $b_{l}(x)<b_{l'}(x)$ (or $v_{l,f}(x)<v_{l',f}(x)$), $\forall l<l'$. Then, we have
   \begin{equation}\label{eq:21}
    v_{1}+\sum_{l=1}^{m}\rho_{l}\frac{\partial v_{l}}{\partial\rho}<\lambda_{1}<v_{1}<\lambda_{2}<...<v_{l-1}<\lambda_{l}<...<v_{m-1}<\lambda_{m}<v_{m},
    \;\text{for}\;u/a\in D,
     \end{equation}
where $D=\{u/a|~\rho_l>0, \forall l, \rho<\rho_{jam}\}$, and $\{\lambda_l\}_{l=1}^m$ are $m$ distinct eigenvalues of system (\ref{eq:20}) for $\theta$ being fixed. However, these eigenvalues cannot be explicitly solved for $m>4$. See \cite{Zhang:2006,Zhang:2008,Donat:2008} for detailed discussions.

To implement the scheme in Section 2, we verify that Eqs. (\ref{eq:10})-(\ref{eq:12}) are applicable. We denote by $\delta_{j+1/2}\rho^-_{l,j+1/2}=\delta_{j+1/2}u^-_{l,j+1/2}/a_{j+1/2}$, $\delta_{j+1/2}\rho^-_{j+1/2}=\sum_{l=1}^m\delta_{j+1/2}\rho^-_{l,j+1/2}$, and $\alpha_{l,j}=(b_{l,j}a_j)/(b_{l,j+1/2}a_{j+1/2})$. For simplicity, we drop the subscript \textquotedblleft $j+1/2$" in the following equations, thus the first equation of (\ref{eq:10}) is equivalent to
   \begin{equation}\label{eq:22}
   \gamma=\frac{\delta\rho^-_{l}v(\delta\rho^-)}
   {\alpha_{l,j}\rho^-_{l}v(\rho^-)},\quad 1\leq l\leq m.
   \end{equation}
By adding up all numerators and denominators over $l$, respectively, we have
 \begin{equation*}
\gamma=\frac{\delta\rho^- v(\delta\rho^-)}{v(\rho^-)\sum_{l=1}^{m}\alpha_{l,j}\rho^-_{l}}\leq\frac{q(\rho^{*})}{v(\rho^-)\sum_{l=1}^{m}\alpha_{l,j}\rho^-_{l}},
\end{equation*}
where $q(\rho^{*})$ is the maximum of the flux function $q(\rho)\equiv\rho v(\rho)$. Thus, the maximum of $\gamma$ reads:
\begin{equation*}
 \gamma_{max}=\min(1, \frac{q(\rho^{*})}{v(\rho)\sum_{l=1}^{m}\alpha_{l}\rho_{l}}).
\end{equation*}
By the inequality of (\ref{eq:21}), the eigenvalues $\{\lambda_l\}_{l=2}^m$ are non-negative, and $\lambda_1$ changes sign such that $\lambda_1\gtrless0$, if $\rho\lessgtr\rho^*$. See \cite{Zhang:2006,Zhang:2008} for detailed discussions. By setting $\gamma=\gamma_{max}$, this together with Eq. (\ref{eq:22}) helps uniquely determine $\delta_{j+1/2}\rho^-_{l,j+1/2}$, under the restrains of Eq. (\ref{eq:11})-(\ref{eq:12}).  Note that $\delta_{j+1/2}\rho^+_{l,j+1/2}$ in Eqs. (\ref{eq:10})-(\ref{eq:12}) and the mapped values in Eqs. (\ref{eq:14})-(\ref{eq:17}) can be similarly derived.

The scheme is applied to resolve complex wave breaking of the Riemann solution, which was analytically discussed and numerically solved in \cite{Zhang:2008}. The initial or Riemann data are given as follows:
\begin{align}\label{eq:23}
u(x,0)=\left \{\aligned
 &u^{L},\;\;\text{if}\;x<x_{0},\\
 &u^{R},\;\;\text{if}\;x>x_{0},
 \endaligned
 \right. \;\; \theta(x)=\left \{\aligned
 &\theta^{L},\;\;\text{if}\;x<x_{0},\\
 &\theta^{R},\;\;\text{if}\;x>x_{0}.
 \endaligned
 \right.
\end{align}
For comparison, we set the same values of all model parameters as those in \cite{Zhang:2008}, and the following values
   \begin{eqnarray*}
   m=3, \;v_{f}=40m/s,\; L=10000m,\; N\equiv L/\triangle_j=800, \;T=400s; \\
   (b^{L}_{1},b^{L}_{2},b^{L}_{3})=(0.5,0.75,1), \quad(b^{R}_{1},b^{R}_{2},b^{R}_{3})=(0.25,0.375,0.5),
   \end{eqnarray*}
are commonly applicable. Others such as $\rho^L\equiv u^L/a^L$, $\rho^R\equiv u^R/a^R$, $r\equiv a^L/a^R$, and $x_0$ in Eq. (\ref{eq:23}) are given in the caption. The eigenvalues used in Eq. (\ref{eq:17}) are estimated by Eq. (\ref{eq:21}). However, we choose the Courant number $C=0.3$, which is smaller than that in \cite{Zhang:2008}. This is under the restriction of Eq. (\ref{eq:17}).

\begin{center}
\scriptsize(a)\epsfig{figure=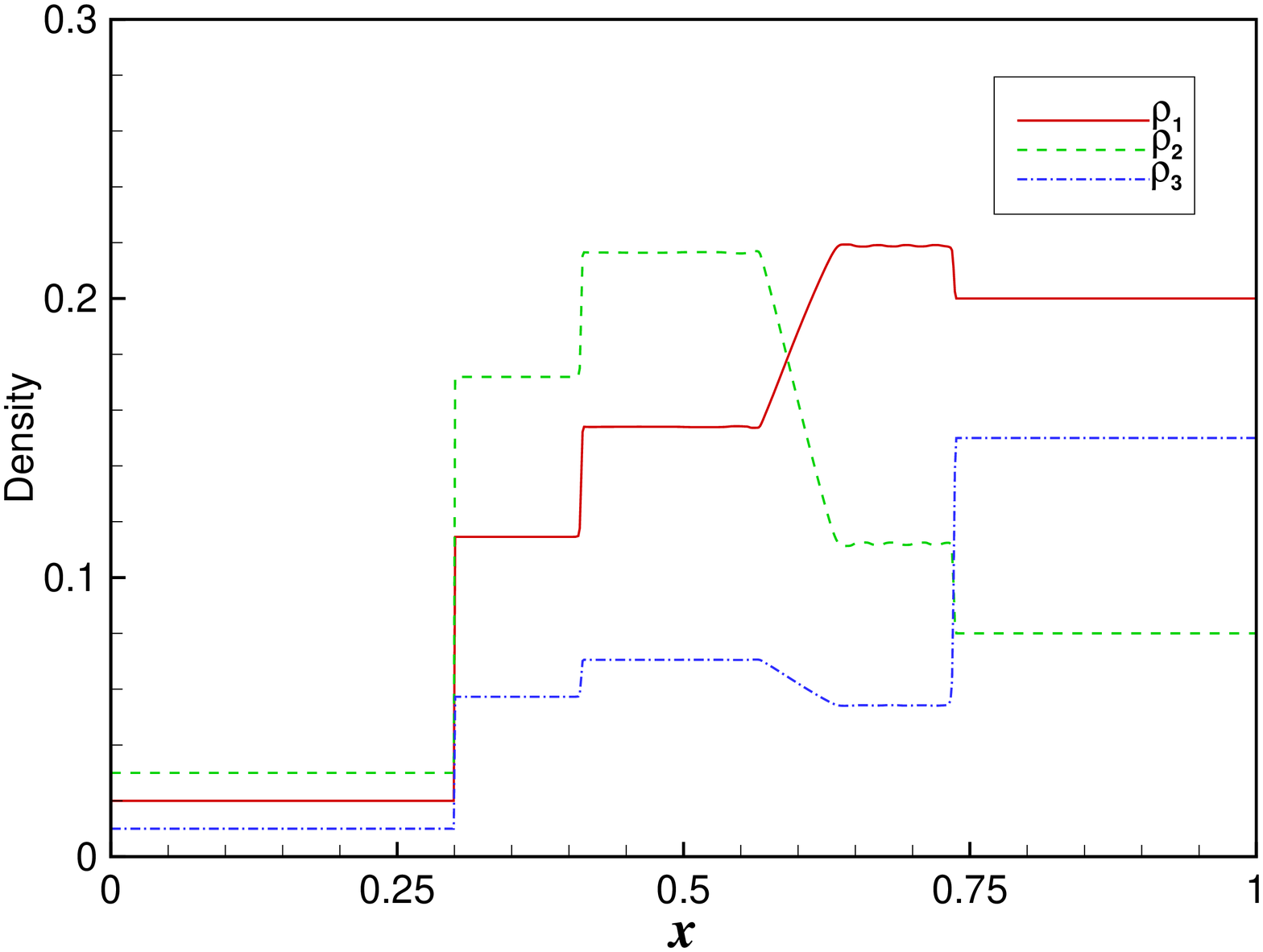,width=3.0 in}  \hspace{0.0 cm}
\scriptsize(b)\epsfig{figure=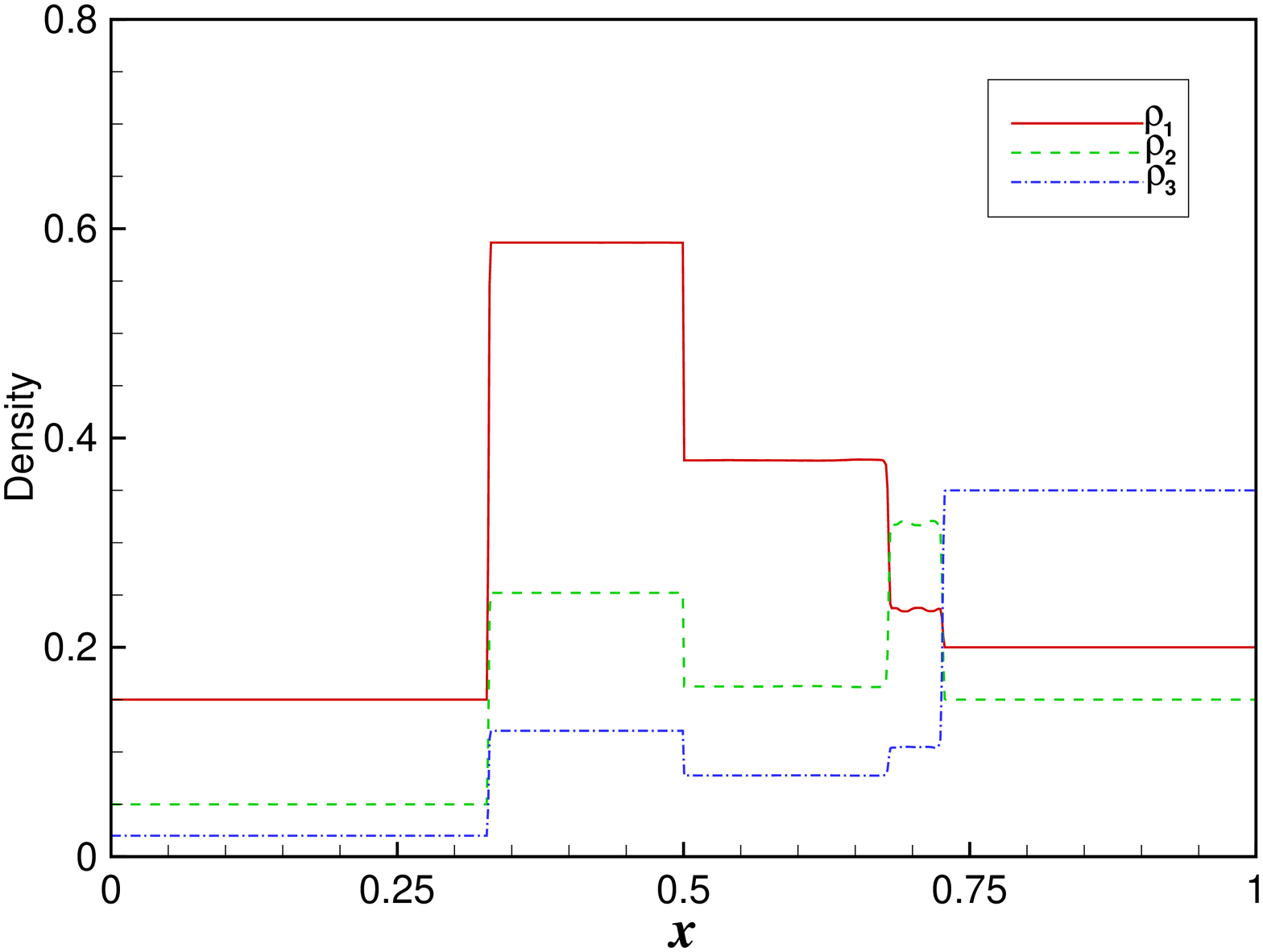,width=3.0 in}  \hspace{0.0 cm}\\
\vspace{0.0cm}
 {\footnotesize Fig. 4. Two wave breaking patterns associated with being strictly hyperbolic of system (\ref{eq:20}), with $\{\rho_{l}\}_{l=1}^3$ being shown at $t=400s$, and (a) $x_{0}=0.3$, $r=2$, $u^{L}/a^{L}=(0.02,0.03,0.01)^{T}$, $u^{R}/a^{R}=(0.2,0.08,0.15)^{T}$, and $\bar{\theta}(\theta_j,\theta_{j+1})=\theta_{j+1}$; (b) $x_{0}=0.5$, $r=3$, $u^{L}/a^{L}=(0.15,0.05,0.02)^{T}$, $u^{R}/a^{R}=(0.2,0.15,0.35)^{T}$, and $\bar{\theta}(\theta_j,\theta_{j+1})=\theta_{j}$.}
\end{center}

\begin{center}
\scriptsize(a)\epsfig{figure=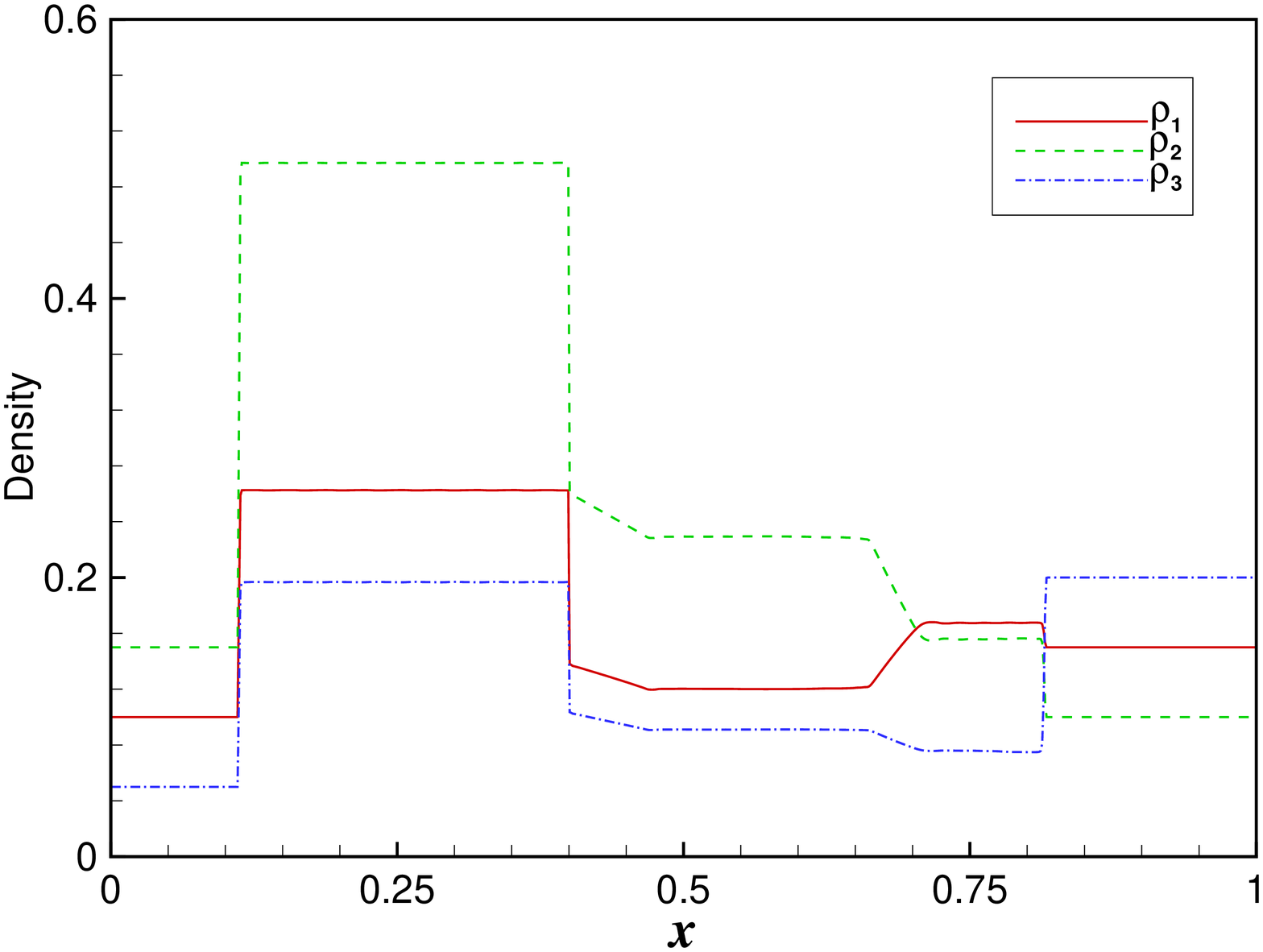,width=3.0 in}  \hspace{0.0 cm}
\scriptsize(b)\epsfig{figure=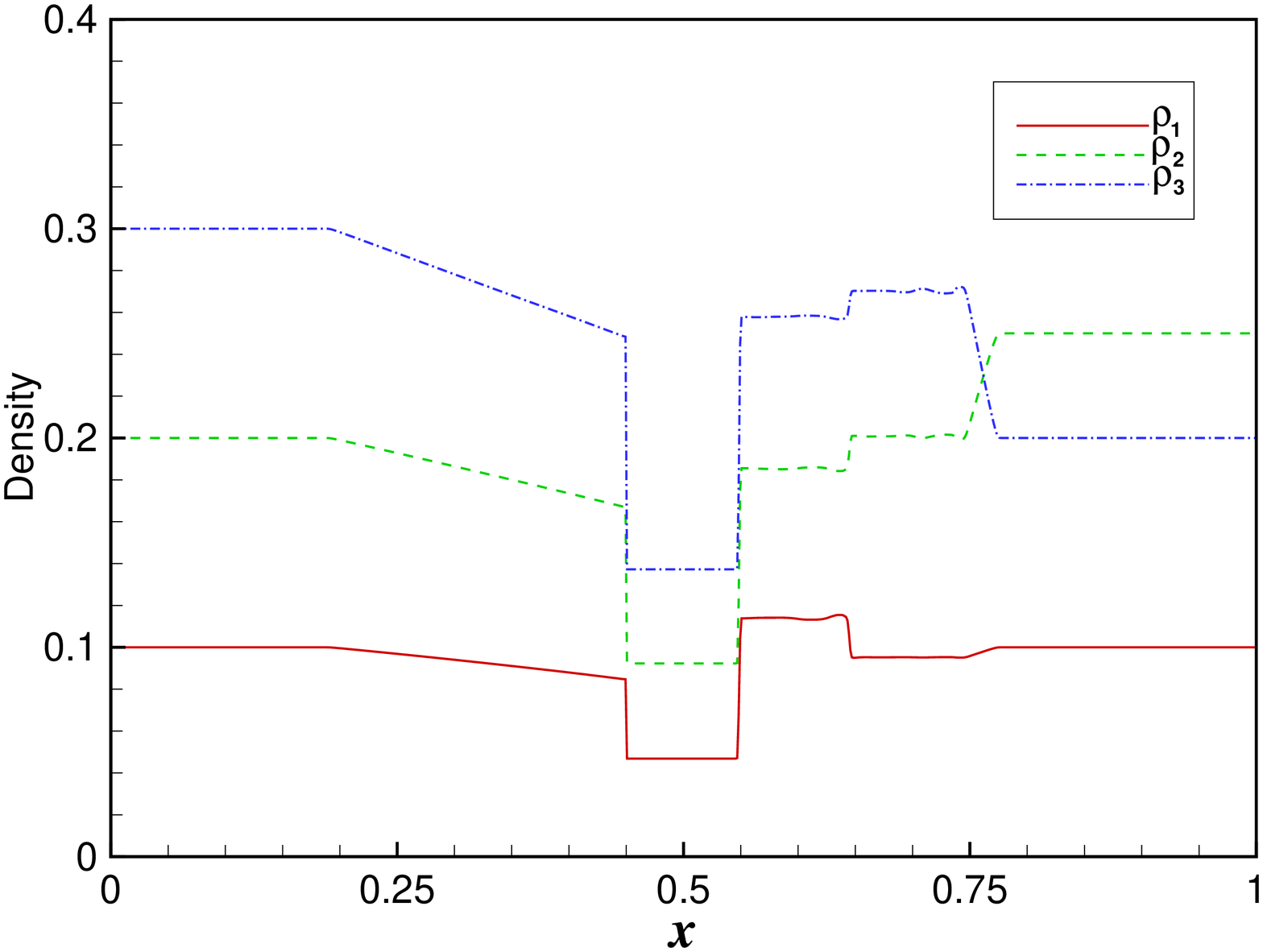,width=3.0 in}  \hspace{0.0 cm}\\
 \vspace{0.0cm}
 {\footnotesize Fig. 5. Two wave breaking patterns associated with being non-strictly hyperbolic of system (\ref{eq:20}), with $\{\rho_{l}\}_{l=1}^3$ being shown at $t=400s$, and (a) $x_{0}=0.4$, $r=3$, $u^{L}/a^{L}=(0.1,0.15,0.05)^{T}$, $u^{R}/a^{R}=(0.15,0.1,0.2)^{T}$, and $\bar{\theta}(\theta_j,\theta_{j+1})=\theta_{j}$; (b) $x_{0}=0.45$, $r=0.4$, $u^{L}/a^{L}=(0.1,0.2,0.3)^{T}$, $u^{R}/a^{R}=(0.1,0.25,0.2)^{T}$, and $\bar{\theta}(\theta_j,\theta_{j+1})=\theta_{j}$.}
\end{center}

%\begin{center}
%\epsfig{figure=6a.eps,width=3.2 in}  \hspace{0.0 cm}
%\epsfig{figure=6b.eps,width=3.2 in}  \hspace{0.0 cm}\\
% \vspace{0.3cm}
% {\footnotesize Fig. 6. The wave breaking pattern of the Riemann problem at $t=400s$,: (a) ; (b) total density $\rho$.}
%\end{center}
%
%\begin{center}
%\epsfig{figure=7a.eps,width=3.2 in}  \hspace{0.0 cm}
%\epsfig{figure=7b.eps,width=3.2 in}  \hspace{0.0 cm}\\
% \vspace{0.3cm}
% {\footnotesize Fig. 7. The wave breaking pattern of the Riemann problem at $t=400s$,: (a) ; (b) total density $\rho$.}
%\end{center}

Figure 4 shows two wave breaking patterns corresponding to a strictly hyperbolic system of Eq. (\ref{eq:20}). As normally we observe $m+1=4$ waves: besides $m$ waves that are associated with $\lambda_l$-characteristic fields, for $l=1,2,3$, there is a contact at $x=x_0$. The inequalities of (\ref{eq:21}) suggest that the characteristic speeds $\{\lambda_l\}_{l=2}^m$ are always non-negative, thus the related waves (shocks or rarefaction fans) arise in the downstream of $x=x_0$. However, $\lambda_1$ can be positive or negative. In Fig. 4(a), we see that the $\lambda_1$-characteristics are all positive, thus those emitting from $t=0$, for $x<x_0\equiv0.3$, are able to pass through the contact $x=x_0$. In contrast, we see that the $\lambda_1$-characteristics in Fig. 4(b) are all negative, thus those emitting from $t=0$, for $x>x_0\equiv0.5$, are able to pass through the contact $x=x_0$.

Figure 5 shows two wave breaking patterns corresponding to a non-strictly hyperbolic system of Eq. (\ref{eq:2}). We observe $m+2$ waves in Fig. 5. In this case, $\lambda_1$-characteristic field \textquotedblleft abnormally" suggests two waves, which propagate downstream and upstream from the contact $x=x_0$, respectively. In Fig. 5(a), the $\lambda_1$-characteristics from the left side of $x=x_0$ are not able to pass through $x=x_0\equiv0.4$, which (with $\lambda_1<0$, for $x=x_0^-$) triggers a return $\lambda_1$-wave on the upstream road. This meanwhile suggests that $\lambda_1=0$, for $x=x_0^+$, which gives rise to a $\lambda_1$-rarefaction on the downstream road. In Fig. 5(b), it is symmetric that the $\lambda_1$-characteristics from the right side are not able to pass through $x=x_0\equiv0.45$, which suggests $\lambda_1>0$, for $x=x_0^+$, and $\lambda_1=0$, for $x=x_0^-$. As a consequence, we observe a return $\lambda_1$-wave and a $\lambda_1$-rarefaction on the downstream and upstream roads, respectively.

See \cite{Zhang:2006,Zhang:2008} for more details about the wave breaking properties of the system. Figs. 4 and 5 together with sufficiently more numerical tests indicate that the proposed hybrid scheme is able to resolve these complex waves efficiently.

\section{Conclusions}
Hyperbolic conservation laws with discontinuous fluxes involve many applications. In the precious works \cite{Zhang:2003,Zhang:2005a,Zhang:2005b}, we proceeded a thorough study of characteristics for the scalar case of Eq. (\ref{eq:2}). For this equation, it is the flow $f(u,\theta)$ (other than the solution variable $u$) remains constant in a characteristic. Therefore, a $\delta$-mapping algorithm was proposed to map each of the involved cell values $u_i$ onto the cell interface by trying to equalize the two flows at $\theta_j$ and $\theta_{j+1/2}$. For a system of Eq. (\ref{eq:2}), the $\delta$-mapping also tries to equalize the two flows, since the propagation of all characteristics (or the \textquotedblleft upwinding") should locally resolve a steady-steady solution or stationary shock, although usually the characteristic fields of the system cannot be analytically solved. In the case that the two flows cannot be equalized, the mapped value $\delta_{j+1/2}u_i$ is determined such that the flow on the cell interface is \textquotedblleft somehow" maximized \cite{Zhang:2008}.

By the mapping, we say that the solution values are \textquotedblleft unified" on the cell interface, and that the system is locally \textquotedblleft standardized". Thus, any classical schemes for Eq. (\ref{eq:1}) can be used to solve Eq. (\ref{eq:2}) just by replacing all involved values $u_i$ with their mapped values $\delta_{j+1/2}u_i$ in the scheme. The present paper enhances this belief in that the adopted RKDG scheme highly resolves complex waves in two application problems, which seems to conclude our studies in this trend. However, the so-called $\delta$-mapping or the underlying concept for upwinging should find more application problems for further validation or improvement, especially when the characteristics cannot reach the cell interface and the flow for determining a mapped value would be \textquotedblleft somehow" maximized.

More relevantly, the maximization would inevitably occur in any traffic flow models (e.g., in higher-order models \cite{Kerner:1994,Lebacque:2007,Siebel:2006,Tang:2008,Bogdanova:2015}) by associating these models with Eq. (\ref{eq:2}), because their first characteristic speed must be allowed to change from negative to positive or vice versa, so as to reflect the dissipation or formation of a traffic jam. By considering traffic flow on a road network (e.g., see \cite{Coclite:2005,Lebacque:1996,Dong:2015,Lin:2015}), the maximization would take place at a junction for finding the \textquotedblleft mapped values" of the neighboring cell values on all incoming and outgoing roads. In this regard, multi-class traffic flow on a road network poses a challenging for the future study.

\bigskip

\noindent \Large \textbf{Acknowledgements} \normalsize

The study was jointly supported by grants from the National Natural Science Foundation of China (11272199), the National Basic Research Program of China (2012CB725404), the Shanghai Program for Innovative Research Team in Universities, the Research Grants Council of the Hong Kong Special Administrative
Region, China (Project No. 17208614) and a National Research Foundation of Korea grant funded by the Korean government (MEST) (NRF-2010-0029446).

\bigskip
{\small \setlength{\baselineskip}{10pt} \setlength{\parskip}{2pt
plus1pt
  minus2pt}


\begin{thebibliography}{999}

\bibitem{Bale:2003}
D.S. Bale, R.J. LeVeque, S. Mitran, J.A. Rossmanith, A wave propagation method for conservation laws and balance laws with spatially varying flux functions, SIAM J. Sci. Comput. 24(3)(2003) 955-978.

\bibitem{Bassi:1997}
F. Bassi, S. Rebay, A high-order accurate discontinuous finite element method for the numerical solution of  the compressible Navier每Stokes equations, J. Comput. Phys. 131(1997) 267-279.

\bibitem{Bogdanova:2015}
A. Bogdanova, M.N. Smirnova, Z.J. Zhu, N.N. Smirnov, Exploring peculiarities of traffic flows with a viscoelastic model, Transportmetrica A: Transport Science 11(7)(2015) 561-578.

\bibitem{Basson:2009}
D.K. Basson, S. Berres, R. B邦rger, On models of polydisperse sedimentation with particle-size-specific hindered-settling factors, Appl. Math. Model. 33(4)(2009) 1815-1835.

\bibitem{Burger:2010a}
R. B\"{u}rger, K.H. Karlsen, J.D. Towers, On some difference schemes and entropy conditions for a class of multi-species kinematic flow models with discontinuous flux, Netw. Heterog. Media 5(3)(2010) 461-485.

\bibitem{Burger:2010b}	
R. B{\"u}rger, R. Donat, P. Mulet, C.A. Vega, Hyperbolicity analysis of polydisperse sedimentation models via a secular equation for the flux Jacobian,
SIAM J. Appl. Math. 70(7)(2010) 2186-2213.

\bibitem{Burger:2011}
R. B{\"u}rger, R. Donat, P. Mulet, C.A. Vega, On the implementation of WENO schemes for a class of polydisperse sedimentation models, J. Comput. Phys. 230(6)(2011) 2322-2344.

\bibitem{Burger:2013}
R. B\"{u}rger, P. Mulet, L. M. Villada, A diffusively corrected multiclass Lighthill-Whitham-Richards traffic model with anticipation lengths and reaction times, Adv. Appl. Math. Mech. 5(5)(2013) 728-758.

\bibitem{Chen:2012}
J.Z. Chen, Z.K. Shi, Y.M. Hu, Numerical solutions of a multiclass traffic flow model on an inhomogeneous highway using a high-resolution relaxed scheme, J. Zhejiang Univ. Sci. C (Computers and Electronics). 13(2012) 29-36.

\bibitem{Cockburn:2001}
B. Cockburn, C.W. Shu, Runge-Kutta discontinuous Galerkin methods for convection-dominated problems, J. Sci. Comput. 16(3)(2001) 173-261.

\bibitem{Coclite:2005}
G.M. Coclite, M. Garavello, B. Piccoli, Traffic flow on a road network, SIAM J. Math. Anal. 36(2005) 1862-1886.

\bibitem{Donat:2008}
R. Donat, P. Mulet, Characteristic-based schemes for multi-class Lighthill-Whitham-Richards traffic models, J. Sci. Comput. 37(3)(2008) 233-250.

\bibitem{Donat:2010}
R. Donat, P. Mulet, A secular equation for the Jacobian matrix of certain multispecies kinematic flow models, Numer. Methods Part. Differ. Eqs. 26(1)(2010) 159-175.

\bibitem{Dong:2015}
C. Dong, Z. Xiong, C. Shao, H. Zhang, A spatialtemporal-based state space approach for freeway network traffic flow modelling and prediction, Transportmetrica A: Transport Science 11(7)(2015) 547-560.

\bibitem{Gottlieb:2001}
S. Gottlieb, C.W. Shu, E. Tadmor, Strong stability-preserving high-order time discretization methods, SIAM Rev. 43(1)(2001) 89-112.

\bibitem{Jin:2003}
W.L. Jin, H.M. Zhang, On the distribution schemes for determining flows through a merge,
Trans. Res. Part B 37(6)(2003) 521-540.

\bibitem{Jin:2009}
W.L. Jin, L. Chen, E.G. Puckett, Supply-demand diagrams and a new framework for analyzing the inhomogeneous Lighthill-Whitham-Richards model, Transportation and Traffic Theory: Golden Jubilee (2009) 603-635.

\bibitem{Kerner:1994}
B.S. Kerner, P. Konhauser, Structure and parameters of clusters
in traffic flow, Phys. Rev. E. 50(1994) 54-83.

\bibitem{Van:2013}
F. van Wageningen-Kessels, B. van't Hof, S.P. Hoogendoorn, H. van Lint, K. Vuik, Anisotropy in generic multi-class traffic flow models, Transportmetrica A: Transport Science 9(5)(2013) 451每472.

\bibitem{Kubatko:2007}
E.J. Kubatko, J.J. Westerink, C. Dawson, Semi discrete discontinuous Galerkin methods and stage-exceeding-order, strong-stability-preserving Runge-Kutta time discretizations, J. Comput. Phys. 222(2007) 832-848.

\bibitem{Lebacque:1996}
J.P. Lebacque, The Godunov scheme and what it means for first order traffic flow models, In: J.B. Lesort (Ed.), Proceedings of the Thirteenth International Symposium on Transportation and Traffic Theory, Lyon, France, (1996) pp. 647-677.

\bibitem{Lebacque:2007}
J.P. Lebacque, S. Mammar, H. Haj-Salem, The Aw-Rascle
and Zhang's model: vacuum problems, existence and regularity of the
solutions of the Riemann problem. Trans. Res. Part B,
41(7)(2007) 710-721.

\bibitem{LeVeque:2002a}
R.J. LeVeque, Finite Volume Methods for Hyperbolic Problems [M]. Cambridge: Cambridge University Press. (2002) 158-200.

\bibitem{LeVeque:2002b}
R.J. LeVeque, Finite-volume methods for non-linear elasticity in heterogeneous media, Internat. J. Numer. Methods Fluids 40(2002) 93-104.

\bibitem{LeVeque:2003}
R.J. LeVeque, D.H. Yong, Solitary waves in layered nonlinear media, SIAM J. Appl. Math. 63(2003) 1539-1560.

\bibitem{Lin:2015}
Z.Y. Lin, P. Zhang, L.Y. Dong, S.C. Wong, K. Choi, Traffic flow on a road network using a conserved higher-order model, In: PROCEEDINGS OF THE INTERNATIONAL CONFERENCE ON NUMERICAL ANALYSIS AND APPLIED MATHEMATICS 2014 (ICNAAM-2014), AIP Publishing, Rhodes, Greece, 1648(2015) pp. 530006.

\bibitem{Ngoduy:2010}
D. Ngoduy, Multiclass first-order modelling of traffic networks using discontinuous flow-density relationships, Transportmetrica 6(2)(2010) 121-141.

\bibitem{Ngoduy:2011}
D. Ngoduy, Multiclass first-order traffic model using stochastic fundamental diagrams, Transportmetrica 7(2)(2011) 111-125.

\bibitem{Qiao:2014}
D.L. Qiao, P. Zhang, S.C. Wong, K. Choi, Discontinuous Galerkin finite element scheme for a conserved higher-order
traffic flow model by exploring Riemann solvers, Appl. Math. Comput. 244(1)(2014) 567-576.

\bibitem{Shu:1998}
C.W. Shu, Essentially non-oscillatory and weighted essentially non-oscillatory schemes for hyperbolic
conservation laws, in: B. Cockburn, C. Johnson, C.W. Shu, E. Tadmor, A. Quarteroni (Eds.), Advanced Numerical Approximation of Nonlinear Hyperbolic Equations, Lecture Notes in Mathematics, Springer, Berlin, 1697(1998) pp. 325-432.

\bibitem{Siebel:2006}
F. Siebel, W. Mauser, On the fundamental
diagram of traffic flow, SIAM J. Appl. Math. 66(2006) 1150-1162.

\bibitem{Tang:2008}
T.Q. Tang, H.J. Huang, G. Xu, A new macro model with consideration of the traffic interruption probability, Phys. A: Statistical Mechanics and its Applications 387(27)(2008) 6845-6856.

\bibitem{Toro:1999}
E.F. Toro, Riemann Solvers and Numerical Methods for Fluid Dynamics [M]. 2nd. Berlin/Heidelberg:
Springer, (1999) 409-520.

\bibitem{Wang:2011}
G. Wang, An Engquist每Osher type finite difference scheme with a discontinuous flux function in space, J. Comput. Appl. Math. 235(17)(2011) 4966-4977.

\bibitem{Wiens:2013}
J.K. Wiens, J.M. Stockie, J.F. Williams, Riemann solver for a kinematic wave traffic model with discontinuous flux, J. Comput. Phys. 242(2013) 1-23.

\bibitem{Wong:2002}
G.C.K. Wong, S.C. Wong, A multi-class traffic flow model-an extension of LWR model with heterogeneous drivers, Trans. Res. Part A 36(2002) 827每841.

\bibitem{Xu:2007}
Z.L. Xu, P. Zhang, R.X. Liu, $\delta$-mapping algorithm coupled with WENO reconstruction for nonlinear elasticity in heterogeneous media. Appl. Numer. Math. 57(1)(2007) 103-116.

\bibitem{ZhangM:2003}
M. Zhang, C.W. Shu, G.C.K. Wong, S.C. Wong, A weighted essentially non-oscillatory numerical scheme for a multi-class Lighthill-Whitham-Richards traffic flow model, J. Compu. Phys. 191(2)(2003) 639-659.

\bibitem{Zhang:2003}
P. Zhang, R.X. Liu, Hyperbolic conservation laws with space-dependent flux: I. Characteristics theory and Riemann problem, J. Comput. Appl. Math. 156(2003) 1-21.

\bibitem{Zhang:2005a}
P. Zhang, R.X. Liu, Hyperbolic conservation laws with space-dependent flux: II. General study of numerical fluxes, J. Comput. Appl. Math. 176(2005) 105-129.

\bibitem{Zhang:2005b}
P. Zhang, R.X. Liu, Generalization of Runge-Kutta discontinuous Galerkin method to LWR traffic flow model with inhomogeneous road conditions, Numer. Methods Part. Differ. Eqs. 21(1)(2005) 80-88.

\bibitem{Zhang:2006}
P. Zhang, R.X. Liu, S.C. Wong, S.Q. Dai. Hyperbolic and kinematic waves of a class of multi-population partial different equations, Euro. J. Appl. Math. 17(2006) 171-200.

\bibitem{Zhang:2008}
P. Zhang, S. C. Wong, Z.L. Xu, A hybrid scheme for solving a multi-class traffic flow model with complex wave breaking. Comput. Methods Appl. Mech. Engrg. 197(45)(2008) 3816-3827.

\end{thebibliography}
\end{document}